\documentclass{amsart}
\usepackage{latexsym}
\usepackage{amssymb}
\usepackage{amsmath}

\newcommand{\tlowername}[2]%
{$\stackrel{\makebox[1pt]{#1}}%
{\begin{picture}(0,0)%
\put(0,0){\makebox(0,6)[t]{\makebox[1pt]{$#2$}}}%
\end{picture}}$}%

\newcommand{\AR}[1]%
{\begin{picture}(#1,0)%
\put(0,0){\vector(1,0){#1}}%
\end{picture}}%

\newcommand{\DOTAR}[1]%
{\NUMBEROFDOTS=#1%
\divide\NUMBEROFDOTS by 3%
\begin{picture}(#1,0)%
\multiput(0,0)(3,0){\NUMBEROFDOTS}{\circle*{1}}%
\put(#1,0){\vector(1,0){0}}%
\end{picture}}%

\newcommand{\MONO}[1]%
{\begin{picture}(#1,0)%
\put(0,0){\vector(1,0){#1}}%
\put(2,-2){\line(0,1){4}}%
\end{picture}}%

\newcommand{\EPI}[1]%
{\begin{picture}(#1,0)(-#1,0)%
\put(-#1,0){\vector(1,0){#1}}%
\put(-6,-2){\line(0,1){4}}%
\end{picture}}%

\newcommand{\BIMO}[1]%
{\begin{picture}(#1,0)(-#1,0)%
\put(-#1,0){\vector(1,0){#1}}%
\put(-6,-2){\line(0,1){4}}%
\put(-#1,-2){\hspace{2pt}\line(0,1){4}}%
\end{picture}}%

\newcommand{\BIAR}[1]%
{\begin{picture}(#1,4)%
\put(0,0){\vector(1,0){#1}}%
\put(0,4){\vector(1,0){#1}}%
\end{picture}}%

\newcommand{\EQL}[1]%
{\begin{picture}(#1,0)%
\put(0,1){\line(1,0){#1}}%
\put(0,-1){\line(1,0){#1}}%
\end{picture}}%

\newcommand{\ADJAR}[1]%
{\begin{picture}(#1,4)%
\put(0,0){\vector(1,0){#1}}%
\put(#1,4){\vector(-1,0){#1}}%
\end{picture}}%

\newcommand{\BKAR}[1]%
{\begin{picture}(#1,0)%
\put(#1,0){\vector(-1,0){#1}}%
\end{picture}}%

\newcommand{\BKDOTAR}[1]%
{\NUMBEROFDOTS=#1%
\divide\NUMBEROFDOTS by 3%
\begin{picture}(#1,0)%
\multiput(#1,0)(-3,0){\NUMBEROFDOTS}{\circle*{1}}%
\put(0,0){\vector(-1,0){0}}%
\end{picture}}%

\newcommand{\BKMONO}[1]%
{\begin{picture}(#1,0)(-#1,0)%
\put(0,0){\vector(-1,0){#1}}%
\put(-2,-2){\line(0,1){4}}%
\end{picture}}%

\newcommand{\BKEPI}[1]%
{\begin{picture}(#1,0)%
\put(#1,0){\vector(-1,0){#1}}%
\put(6,-2){\line(0,1){4}}%
\end{picture}}%

\newcommand{\BKBIMO}[1]%
{\begin{picture}(#1,0)%
\put(#1,0){\vector(-1,0){#1}}%
\put(6,-2){\line(0,1){4}}%
\put(#1,-2){\hspace{-2pt}\line(0,1){4}}%
\end{picture}}%

\newcommand{\BKBIAR}[1]%
{\begin{picture}(#1,4)%
\put(#1,0){\vector(-1,0){#1}}%
\put(#1,4){\vector(-1,0){#1}}%
\end{picture}}%

\newcommand{\BKADJAR}[1]%
{\begin{picture}(#1,4)%
\put(0,4){\vector(1,0){#1}}%
\put(#1,0){\vector(-1,0){#1}}%
\end{picture}}%

\newcommand{\lowername}[2]%
{$\stackrel{\makebox[1pt]{#1}}%
{\begin{picture}(0,0)%
\truex{600}%
\put(0,0){\makebox(0,\value{x})[t]{\makebox[1pt]{$#2$}}}%
\end{picture}}$}%

\newcommand{\hcase}[2]%
{\makebox[0pt]%
{\raisebox{-1pt}[0pt][0pt]{#1{#2}}}}%

\newcommand{\Hcase}[3]%
{\makebox[0pt]
{\raisebox{-1pt}[0pt][0pt]%
{$\stackrel{\makebox[0pt]{$\textstyle{#2}$}}{#1{#3}}$}}}%

\newcommand{\hcasE}[3]%
{\makebox[0pt]%
{\raisebox{-9pt}[0pt][0pt]%
{\lowername{#1{#3}}{#2}}}}%

\newcommand{\hbicase}[2]%
{\makebox[0pt]%
{\raisebox{-2.5pt}[0pt][0pt]{#1{#2}}}}%

\newcommand{\Hbicase}[4]%
{\makebox[0pt]
{\raisebox{-10.5pt}[0pt][0pt]%
{$\stackrel{\makebox[0pt]{$\textstyle{#2}$}}%
{\mbox{\lowername{#1{#4}}{#3}}}$}}}%

\newcommand{\EAR}[1]%
{\begin{picture}(#1,0)%
\put(0,0){\vector(1,0){#1}}%
\end{picture}}%

\newcommand{\EDOTAR}[1]%
{\truex{100}\truey{300}%
\NUMBEROFDOTS=#1%
\divide\NUMBEROFDOTS by \value{y}%
\begin{picture}(#1,0)%
\multiput(0,0)(\value{y},0){\NUMBEROFDOTS}%
{\circle*{\value{x}}}%
\put(#1,0){\vector(1,0){0}}%
\end{picture}}%

\newcommand{\EMONO}[1]%
{\begin{picture}(#1,0)%
\put(0,0){\vector(1,0){#1}}%
\truex{300}\truey{600}%
\put(\value{x},-\value{x}){\line(0,1){\value{y}}}%
\end{picture}}%

\newcommand{\EEPI}[1]%
{\begin{picture}(#1,0)(-#1,0)%
\put(-#1,0){\vector(1,0){#1}}%
\truex{300}\truey{600}\truez{800}%
\put(-\value{z},-\value{x}){\line(0,1){\value{y}}}%
\end{picture}}%

\newcommand{\EBIMO}[1]%
{\begin{picture}(#1,0)(-#1,0)%
\put(-#1,0){\vector(1,0){#1}}%
\truex{300}\truey{600}\truez{800}%
\put(-\value{z},-\value{x}){\line(0,1){\value{y}}}%
\put(-#1,-\value{x}){\hspace{3pt}\line(0,1){\value{y}}}%
\end{picture}}%

\newcommand{\EBIAR}[1]%
{\truex{400}%
\begin{picture}(#1,\value{x})%
\put(0,0){\vector(1,0){#1}}%
\put(0,\value{x}){\vector(1,0){#1}}%
\end{picture}}%

\newcommand{\EEQL}[1]%
{\begin{picture}(#1,0)%
\truex{200}%
\put(0,\value{x}){\line(1,0){#1}}%
\put(0,0){\line(1,0){#1}}%
\end{picture}}%

\newcommand{\EADJAR}[1]%
{\truex{400}%
\begin{picture}(#1,\value{x})%
\put(0,0){\vector(1,0){#1}}%
\put(#1,\value{x}){\vector(-1,0){#1}}%
\end{picture}}%

\newcommand{\ear}%
{\hspace{\SOURCE\unitlength}%
\hcase{\EAR}{\ARROWLENGTH}}%

\newcommand{\Ear}[1]%
{\hspace{\SOURCE\unitlength}%
\Hcase{\EAR}{#1}{\ARROWLENGTH}}%

\newcommand{\eaR}[1]%
{\hspace{\SOURCE\unitlength}%
\hcasE{\EAR}{#1}{\ARROWLENGTH}}%

\newcommand{\edotar}%
{\hspace{\SOURCE\unitlength}%
\hcase{\EDOTAR}{\ARROWLENGTH}}%

\newcommand{\Edotar}[1]%
{\hspace{\SOURCE\unitlength}%
\Hcase{\EDOTAR}{#1}{\ARROWLENGTH}}%

\newcommand{\edotaR}[1]%
{\hspace{\SOURCE\unitlength}%
\hcasE{\EDOTAR}{#1}{\ARROWLENGTH}}%

\newcommand{\emono}%
{\hspace{\SOURCE\unitlength}%
\hcase{\EMONO}{\ARROWLENGTH}}%

\newcommand{\Emono}[1]%
{\hspace{\SOURCE\unitlength}%
\Hcase{\EMONO}{#1}{\ARROWLENGTH}}%

\newcommand{\emonO}[1]%
{\hspace{\SOURCE\unitlength}%
\hcasE{\EMONO}{#1}{\ARROWLENGTH}}%

\newcommand{\eepi}%
{\hspace{\SOURCE\unitlength}%
\hcase{\EEPI}{\ARROWLENGTH}}%

\newcommand{\Eepi}[1]%
{\hspace{\SOURCE\unitlength}%
\Hcase{\EEPI}{#1}{\ARROWLENGTH}}%

\newcommand{\eepI}[1]%
{\hspace{\SOURCE\unitlength}%
\hcasE{\EEPI}{#1}{\ARROWLENGTH}}%

\newcommand{\ebimo}%
{\hspace{\SOURCE\unitlength}%
\hcase{\EBIMO}{\ARROWLENGTH}}%

\newcommand{\Ebimo}[1]%
{\hspace{\SOURCE\unitlength}%
\Hcase{\EBIMO}{#1}{\ARROWLENGTH}}%

\newcommand{\ebimO}[1]%
{\hspace{\SOURCE\unitlength}%
\hcasE{\EBIMO}{#1}{\ARROWLENGTH}}%

\newcommand{\eiso}%
{\hspace{\SOURCE\unitlength}%
\Hcase{\EAR}{\cong}{\ARROWLENGTH}}%

\newcommand{\Eiso}[1]%
{\hspace{\SOURCE\unitlength}%
\Hcase{\EAR}{\cong#1}{\ARROWLENGTH}}%

\newcommand{\eisO}[1]%
{\hspace{\SOURCE\unitlength}%
\hcasE{\EAR}{\cong#1}{\ARROWLENGTH}}%

\newcommand{\ebiar}%
{\hspace{\SOURCE\unitlength}%
\hbicase{\EBIAR}{\ARROWLENGTH}}%

\newcommand{\Ebiar}[2]%
{\hspace{\SOURCE\unitlength}%
\Hbicase{\EBIAR}{#1}{#2}{\ARROWLENGTH}}%

\newcommand{\eeql}%
{\hspace{\SOURCE\unitlength}%
\hbicase{\EEQL}{\ARROWLENGTH}}%

\newcommand{\eadjar}%
{\hspace{\SOURCE\unitlength}%
\hbicase{\EADJAR}{\ARROWLENGTH}}%

\newcommand{\Eadjar}[2]%
{\hspace{\SOURCE\unitlength}%
\Hbicase{\EADJAR}{#1}{#2}{\ARROWLENGTH}}%

\newcommand{\WAR}[1]%
{\begin{picture}(#1,0)%
\put(#1,0){\vector(-1,0){#1}}%
\end{picture}}%

\newcommand{\WDOTAR}[1]%
{\truex{100}\truey{300}%
\NUMBEROFDOTS=#1%
\divide\NUMBEROFDOTS by \value{y}%
\begin{picture}(#1,0)%
\multiput(#1,0)(-\value{y},0){\NUMBEROFDOTS}%
{\circle*{\value{x}}}%
\put(0,0){\vector(-1,0){0}}%
\end{picture}}%

\newcommand{\WMONO}[1]%
{\begin{picture}(#1,0)(-#1,0)%
\put(0,0){\vector(-1,0){#1}}%
\truex{300}\truey{600}%
\put(-\value{x},-\value{x}){\line(0,1){\value{y}}}%
\end{picture}}%

\newcommand{\WEPI}[1]%
{\begin{picture}(#1,0)%
\put(#1,0){\vector(-1,0){#1}}%
\truex{300}\truey{600}\truez{800}%
\put(\value{z},-\value{x}){\line(0,1){\value{y}}}%
\end{picture}}%

\newcommand{\WBIMO}[1]%
{\begin{picture}(#1,0)%
\put(#1,0){\vector(-1,0){#1}}%
\truex{300}\truey{600}\truez{800}%
\put(\value{z},-\value{x}){\line(0,1){\value{y}}}%
\put(#1,-\value{x}){\hspace{-3pt}\line(0,1){\value{y}}}%
\end{picture}}%

\newcommand{\WBIAR}[1]%
{\truex{400}%
\begin{picture}(#1,\value{x})%
\put(#1,0){\vector(-1,0){#1}}%
\put(#1,\value{x}){\vector(-1,0){#1}}%
\end{picture}}%

\newcommand{\WADJAR}[1]%
{\truex{400}%
\begin{picture}(#1,\value{x})%
\put(0,\value{x}){\vector(1,0){#1}}%
\put(#1,0){\vector(-1,0){#1}}%
\end{picture}}%

\newcommand{\war}%
{\hspace{\SOURCE\unitlength}%
\hcase{\WAR}{\ARROWLENGTH}}%

\newcommand{\War}[1]%
{\hspace{\SOURCE\unitlength}%
\Hcase{\WAR}{#1}{\ARROWLENGTH}}%

\newcommand{\waR}[1]%
{\hspace{\SOURCE\unitlength}%
\hcasE{\WAR}{#1}{\ARROWLENGTH}}%

\newcommand{\wdotar}%
{\hspace{\SOURCE\unitlength}%
\hcase{\WDOTAR}{\ARROWLENGTH}}%

\newcommand{\Wdotar}[1]%
{\hspace{\SOURCE\unitlength}%
\Hcase{\WDOTAR}{#1}{\ARROWLENGTH}}%

\newcommand{\wdotaR}[1]%
{\hspace{\SOURCE\unitlength}%
\hcasE{\WDOTAR}{#1}{\ARROWLENGTH}}%

\newcommand{\wmono}%
{\hspace{\SOURCE\unitlength}%
\hcase{\WMONO}{\ARROWLENGTH}}%

\newcommand{\Wmono}[1]%
{\hspace{\SOURCE\unitlength}%
\Hcase{\WMONO}{#1}{\ARROWLENGTH}}%

\newcommand{\wmonO}[1]%
{\hspace{\SOURCE\unitlength}%
\hcasE{\WMONO}{#1}{\ARROWLENGTH}}%

\newcommand{\wepi}%
{\hspace{\SOURCE\unitlength}%
\hcase{\WEPI}{\ARROWLENGTH}}%

\newcommand{\Wepi}[1]%
{\hspace{\SOURCE\unitlength}%
\Hcase{\WEPI}{#1}{\ARROWLENGTH}}%

\newcommand{\wepI}[1]%
{\hspace{\SOURCE\unitlength}%
\hcasE{\WEPI}{#1}{\ARROWLENGTH}}%

\newcommand{\wbimo}%
{\hspace{\SOURCE\unitlength}%
\hcase{\WBIMO}{\ARROWLENGTH}}%

\newcommand{\Wbimo}[1]%
{\hspace{\SOURCE\unitlength}%
\Hcase{\WBIMO}{#1}{\ARROWLENGTH}}%

\newcommand{\wbimO}[1]%
{\hspace{\SOURCE\unitlength}%
\hcasE{\WBIMO}{#1}{\ARROWLENGTH}}%

\newcommand{\wiso}%
{\hspace{\SOURCE\unitlength}%
\Hcase{\WAR}{\cong}{\ARROWLENGTH}}%

\newcommand{\Wiso}[1]%
{\hspace{\SOURCE\unitlength}%
\Hcase{\WAR}{#1}{\ARROWLENGTH}}%

\newcommand{\wisO}[1]%
{\hspace{\SOURCE\unitlength}%
\hcasE{\WAR}{#1}{\ARROWLENGTH}}%

\newcommand{\wbiar}%
{\hspace{\SOURCE\unitlength}%
\hbicase{\WBIAR}{\ARROWLENGTH}}%

\newcommand{\Wbiar}[2]%
{\hspace{\SOURCE\unitlength}%
\Hbicase{\WBIAR}{#1}{#2}{\ARROWLENGTH}}%

\newcommand{\weql}%
{\hspace{\SOURCE\unitlength}%
\hbicase{\EEQL}{\ARROWLENGTH}}%

\newcommand{\wadjar}%
{\hspace{\SOURCE\unitlength}%
\hbicase{\WADJAR}{\ARROWLENGTH}}%

\newcommand{\Wadjar}[2]%
{\hspace{\SOURCE\unitlength}%
\Hbicase{\WADJAR}{#1}{#2}{\ARROWLENGTH}}%

\newcommand{\vcase}[2]{#1{#2}}%

\newcommand{\Vcase}[3]{\makebox[0pt]%
{\makebox[0pt][r]{\raisebox{0pt}[0pt][0pt]{${#2}\hspace{2pt}$}}}#1{#3}}%

\newcommand{\vcasE}[3]{\makebox[0pt]%
{#1{#3}\makebox[0pt][l]{\raisebox{0pt}[0pt][0pt]{\hspace{2pt}$#2$}}}}%

\newcommand{\vbicase}[2]{\makebox[0pt]{{#1{#2}}}}%

\newcommand{\Vbicase}[4]{\makebox[0pt]%
{\makebox[0pt][r]{\raisebox{0pt}[0pt][0pt]{$#2$\hspace{4pt}}}#1{#4}%
\makebox[0pt][l]{\raisebox{0pt}[0pt][0pt]{\hspace{5pt}$#3$}}}}%

\newcommand{\SAR}[1]%
{\begin{picture}(0,0)%
\put(0,0){\makebox(0,0)%
{\begin{picture}(0,#1)%
\put(0,#1){\vector(0,-1){#1}}%
\end{picture}}}\end{picture}}%

\newcommand{\SDOTAR}[1]%
{\truex{100}\truey{300}%
\NUMBEROFDOTS=#1%
\divide\NUMBEROFDOTS by \value{y}%
\begin{picture}(0,0)%
\put(0,0){\makebox(0,0)%
{\begin{picture}(0,#1)%
\multiput(0,#1)(0,-\value{y}){\NUMBEROFDOTS}%
{\circle*{\value{x}}}%
\put(0,0){\vector(0,-1){0}}%
\end{picture}}}\end{picture}}%

\newcommand{\SMONO}[1]%
{\begin{picture}(0,0)%
\put(0,0){\makebox(0,0)%
{\begin{picture}(0,#1)%
\put(0,#1){\vector(0,-1){#1}}%
\truex{300}\truey{600}%
\put(0,#1){\begin{picture}(0,0)%
\put(-\value{x},-\value{x}){\line(1,0){\value{y}}}\end{picture}}%
\end{picture}}}\end{picture}}%

\newcommand{\SEPI}[1]%
{\begin{picture}(0,0)%
\put(0,0){\makebox(0,0)%
{\begin{picture}(0,#1)%
\put(0,#1){\vector(0,-1){#1}}%
\truex{300}\truey{600}\truez{800}%
\put(-\value{x},\value{z}){\line(1,0){\value{y}}}%
\end{picture}}}\end{picture}}%

\newcommand{\SBIMO}[1]%
{\begin{picture}(0,0)%
\put(0,0){\makebox(0,0)%
{\begin{picture}(0,#1)%
\put(0,#1){\vector(0,-1){#1}}%
\truex{300}\truey{600}\truez{800}%
\put(0,#1){\begin{picture}(0,0)%
\put(-\value{x},-\value{x}){\line(1,0){\value{y}}}\end{picture}}%
\put(-\value{x},\value{z}){\line(1,0){\value{y}}}%
\end{picture}}}\end{picture}}%

\newcommand{\SBIAR}[1]%
{\begin{picture}(0,0)%
\truex{200}%
\put(0,0){\makebox(0,0)%
{\begin{picture}(0,#1)\put(-\value{x},#1){\vector(0,-1){#1}}%
\put(\value{x},#1){\vector(0,-1){#1}}%
\end{picture}}}\end{picture}}%

\newcommand{\SEQL}[1]%
{\begin{picture}(0,0)%
\truex{100}%
\put(0,0){\makebox(0,0)%
{\begin{picture}(0,#1)\put(-\value{x},#1){\line(0,-1){#1}}%
\put(\value{x},#1){\line(0,-1){#1}}%
\end{picture}}}\end{picture}}%

\newcommand{\sarv}[1]{\vcase{\SAR}{#100}}%

\newcommand{\sar}{\sarv{50}}%

\newcommand{\Sarv}[2]{\Vcase{\SAR}{#1}{#200}}%

\newcommand{\Sar}[1]{\Sarv{#1}{50}}%

\newcommand{\saRv}[2]{\vcasE{\SAR}{#1}{#200}}%

\newcommand{\saR}[1]{\saRv{#1}{50}}%

\newcommand{\Sisov}[2]%
{\Vbicase{\SAR}{#1\hspace{-2pt}}{\hspace{-2pt}\cong}{#200}}%

\newcommand{\seqlv}[1]{\vbicase{\SEQL}{#100}}%

\newcommand{\seql}{\seqlv{50}}%

\newcommand{\NAR}[1]%
{\begin{picture}(0,0)%
\put(0,0){\makebox(0,0)%
{\begin{picture}(0,#1)\put(0,0){\vector(0,1){#1}}%
\end{picture}}}\end{picture}}%

\newcommand{\NDOTAR}[1]%
{\truex{100}\truey{300}%
\NUMBEROFDOTS=#1%
\divide\NUMBEROFDOTS by \value{y}%
\begin{picture}(0,0)%
\put(0,0){\makebox(0,0)%
{\begin{picture}(0,#1)%
\multiput(0,0)(0,\value{y}){\NUMBEROFDOTS}%
{\circle*{\value{x}}}%
\put(0,#1){\vector(0,1){0}}%
\end{picture}}}\end{picture}}%

\newcommand{\NMONO}[1]%
{\begin{picture}(0,0)%
\put(0,0){\makebox(0,0)%
{\begin{picture}(0,#1)%
\put(0,0){\vector(0,1){#1}}%
\truex{300}\truey{600}%
\put(-\value{x},\value{x}){\line(1,0){\value{y}}}%
\end{picture}}}%
\end{picture}}%

\newcommand{\NEPI}[1]%
{\begin{picture}(0,0)%
\put(0,0){\makebox(0,0)%
{\begin{picture}(0,#1)%
\put(0,0){\vector(0,1){#1}}%
\truex{300}\truey{600}\truez{800}%
\put(0,#1){\begin{picture}(0,0)%
\put(-\value{x},-\value{z}){\line(1,0){\value{y}}}\end{picture}}%
\end{picture}}}\end{picture}}%

\newcommand{\NBIMO}[1]%
{\begin{picture}(0,0)%
\put(0,0){\makebox(0,0)%
{\begin{picture}(0,#1)%
\put(0,0){\vector(0,1){#1}}%
\truex{300}\truey{600}\truez{800}%
\put(-\value{x},\value{x}){\line(1,0){\value{y}}}%
\put(0,#1){\begin{picture}(0,0)%
\put(-\value{x},-\value{z}){\line(1,0){\value{y}}}\end{picture}}%
\end{picture}}}\end{picture}}%

\newcommand{\NBIAR}[1]%
{\begin{picture}(0,0)%
\truex{200}%
\put(0,0){\makebox(0,0)%
{\begin{picture}(0,#1)\put(-\value{x},0){\vector(0,1){#1}}%
\put(\value{x},0){\vector(0,1){#1}}%
\end{picture}}}\end{picture}}%

\newcommand{\Nisov}[2]%
{\Vbicase{\NAR}{#1\hspace{-2pt}}{\hspace{-2pt}\cong}{#200}}%

\newcommand{\NEDOTAR}%
{\truex{100}\truey{212}%
\NUMBEROFDOTS=5800%
\divide\NUMBEROFDOTS by \value{y}%
\begin{picture}(0,0)%
\multiput(-2900,-2900)(\value{y},\value{y}){\NUMBEROFDOTS}%
{\circle*{\value{x}}}%
\put(2900,2900){\vector(1,1){0}}%
\end{picture}}%

\newcommand{\SWDOTAR}%
{\truex{100}\truey{212}%
\NUMBEROFDOTS=5800%
\divide\NUMBEROFDOTS by \value{y}%
\begin{picture}(0,0)%
\multiput(2900,2900)(-\value{y},-\value{y}){\NUMBEROFDOTS}%
{\circle*{\value{x}}}%
\put(-2900,-2900){\vector(-1,-1){0}}%
\end{picture}}%

\newcommand{\SEDOTAR}%
{\truex{100}\truey{212}%
\NUMBEROFDOTS=5800%
\divide\NUMBEROFDOTS by \value{y}%
\begin{picture}(0,0)%
\multiput(-2900,2900)(\value{y},-\value{y}){\NUMBEROFDOTS}%
{\circle*{\value{x}}}%
\put(2900,-2900){\vector(1,-1){0}}%
\end{picture}}%

\newcommand{\NWDOTAR}%
{\truex{100}\truey{212}%
\NUMBEROFDOTS=5800%
\divide\NUMBEROFDOTS by \value{y}%
\begin{picture}(0,0)%
\multiput(2900,-2900)(-\value{y},\value{y}){\NUMBEROFDOTS}%
{\circle*{\value{x}}}%
\put(-2900,2900){\vector(-1,1){0}}%
\end{picture}}%

\newcommand{\ENEAR}[2]%
{\makebox[0pt]{\begin{picture}(0,0)%
\put(0,-150){\makebox(0,0){\begin{picture}(0,0)%
\put(-6600,-3300){\vector(2,1){13200}}%
\truex{200}\truey{800}\truez{600}%
\put(-\value{x},\value{x}){\makebox(0,\value{z})[r]{${#1}$}}%
\put(\value{x},-\value{y}){\makebox(0,\value{z})[l]{${#2}$}}%
\end{picture}}}\end{picture}}}%

\newcommand{\ESEAR}[2]%
{\makebox[0pt]{\begin{picture}(0,0)%
\put(0,-150){\makebox(0,0){\begin{picture}(0,0)%
\put(-6600,3300){\vector(2,-1){13200}}%
\truex{200}\truey{800}\truez{600}%
\put(\value{x},\value{x}){\makebox(0,\value{z})[l]{${#1}$}}%
\put(-\value{x},-\value{y}){\makebox(0,\value{z})[r]{${#2}$}}%
\end{picture}}}\end{picture}}}%

\newcommand{\WNWAR}[2]%
{\makebox[0pt]{\begin{picture}(0,0)%
\put(0,-150){\makebox(0,0){\begin{picture}(0,0)%
\put(6600,-3300){\vector(-2,1){13200}}%
\truex{200}\truey{800}\truez{600}%
\put(\value{x},\value{x}){\makebox(0,\value{z})[l]{${#1}$}}%
\put(-\value{x},-\value{y}){\makebox(0,\value{z})[r]{${#2}$}}%
\end{picture}}}\end{picture}}}%

\newcommand{\WSWAR}[2]%
{\makebox[0pt]{\begin{picture}(0,0)%
\put(0,-150){\makebox(0,0){\begin{picture}(0,0)%
\put(6600,3300){\vector(-2,-1){13200}}%
\truex{200}\truey{800}\truez{600}%
\put(-\value{x},\value{x}){\makebox(0,\value{z})[r]{${#1}$}}%
\put(\value{x},-\value{y}){\makebox(0,\value{z})[l]{${#2}$}}%
\end{picture}}}\end{picture}}}%

\newcommand{\NNEAR}[2]%
{\raisebox{-1pt}[0pt][0pt]{\begin{picture}(0,0)%
\put(0,0){\makebox(0,0){\begin{picture}(0,0)%
\put(-3300,-6600){\vector(1,2){6600}}%
\truex{100}\truez{600}%
\put(-\value{x},\value{x}){\makebox(0,\value{z})[r]{${#1}$}}%
\put(\value{x},-\value{z}){\makebox(0,\value{z})[l]{${#2}$}}%
\end{picture}}}\end{picture}}}%

\newcommand{\SSWAR}[2]%
{\raisebox{-1pt}[0pt][0pt]{\begin{picture}(0,0)%
\put(0,0){\makebox(0,0){\begin{picture}(0,0)%
\put(3300,6600){\vector(-1,-2){6600}}%
\truex{100}\truez{600}%
\put(-\value{x},\value{x}){\makebox(0,\value{z})[r]{${#1}$}}%
\put(\value{x},-\value{z}){\makebox(0,\value{z})[l]{${#2}$}}%
\end{picture}}}\end{picture}}}%

\newcommand{\SSEAR}[2]%
{\raisebox{-1pt}[0pt][0pt]{\begin{picture}(0,0)%
\put(0,0){\makebox(0,0){\begin{picture}(0,0)%
\put(-3300,6600){\vector(1,-2){6600}}%
\truex{200}\truez{600}%
\put(\value{x},\value{x}){\makebox(0,\value{z})[l]{${#1}$}}%
\put(-\value{x},-\value{z}){\makebox(0,\value{z})[r]{${#2}$}}%
\end{picture}}}\end{picture}}}%

\newcommand{\NNWAR}[2]%
{\raisebox{-1pt}[0pt][0pt]{\begin{picture}(0,0)%
\put(0,0){\makebox(0,0){\begin{picture}(0,0)%
\put(3300,-6600){\vector(-1,2){6600}}%
\truex{200}\truez{600}%
\put(\value{x},\value{x}){\makebox(0,\value{z})[l]{${#1}$}}%
\put(-\value{x},-\value{z}){\makebox(0,\value{z})[r]{${#2}$}}%
\end{picture}}}\end{picture}}}%

\newcommand{\Necurve}[2]%
{\begin{picture}(0,0)%
\truex{1300}\truey{2000}\truez{200}%
\put(0,\value{x}){\oval(#200,\value{y})[t]}%
\put(0,\value{x}){\makebox(0,0){\begin{picture}(#200,0)%
\put(#200,0){\vector(0,-1){\value{z}}}%
\put(0,0){\line(0,-1){\value{z}}}\end{picture}}}%
\truex{2500}%
\put(0,\value{x}){\makebox(0,0)[b]{${#1}$}}%
\end{picture}}%

\newcommand{\Nwcurve}[2]%
{\begin{picture}(0,0)%
\truex{1300}\truey{2000}\truez{200}%
\put(0,\value{x}){\oval(#200,\value{y})[t]}%
\put(0,\value{x}){\makebox(0,0){\begin{picture}(#200,0)%
\put(#200,0){\line(0,-1){\value{z}}}%
\put(0,0){\vector(0,-1){\value{z}}}\end{picture}}}%
\truex{2500}%
\put(0,\value{x}){\makebox(0,0)[b]{${#1}$}}%
\end{picture}}%

\newcommand{\Securve}[2]%
{\begin{picture}(0,0)%
\truex{1300}\truey{2000}\truez{200}%
\put(0,-\value{x}){\oval(#200,\value{y})[b]}%
\put(0,-\value{x}){\makebox(0,0){\begin{picture}(#200,0)%
\put(#200,0){\vector(0,1){\value{z}}}%
\put(0,0){\line(0,1){\value{z}}}\end{picture}}}%
\truex{2500}%
\put(0,-\value{x}){\makebox(0,0)[t]{${#1}$}}%
\end{picture}}%

\newcommand{\Swcurve}[2]%
{\begin{picture}(0,0)%
\truex{1300}\truey{2000}\truez{200}%
\put(0,-\value{x}){\oval(#200,\value{y})[b]}%
\put(0,-\value{x}){\makebox(0,0){\begin{picture}(#200,0)%
\put(#200,0){\line(0,1){\value{z}}}%
\put(0,0){\vector(0,1){\value{z}}}\end{picture}}}%
\truex{2500}%
\put(0,-\value{x}){\makebox(0,0)[t]{${#1}$}}%
\end{picture}}%

\newcommand{\Escurve}[2]%
{\begin{picture}(0,0)%
\truex{1400}\truey{2000}\truez{200}%
\put(\value{x},0){\oval(\value{y},#200)[r]}%
\put(\value{x},0){\makebox(0,0){\begin{picture}(0,#200)%
\put(0,0){\vector(-1,0){\value{z}}}%
\put(0,#200){\line(-1,0){\value{z}}}\end{picture}}}%
\truex{2500}%
\put(\value{x},0){\makebox(0,0)[l]{${#1}$}}%
\end{picture}}%

\newcommand{\Encurve}[2]%
{\begin{picture}(0,0)%
\truex{1400}\truey{2000}\truez{200}%
\put(\value{x},0){\oval(\value{y},#200)[r]}%
\put(\value{x},0){\makebox(0,0){\begin{picture}(0,#200)%
\put(0,0){\line(-1,0){\value{z}}}%
\put(0,#200){\vector(-1,0){\value{z}}}\end{picture}}}%
\truex{2500}%
\put(\value{x},0){\makebox(0,0)[l]{${#1}$}}%
\end{picture}}%

\newcommand{\Wscurve}[2]%
{\begin{picture}(0,0)%
\truex{1300}\truey{2000}\truez{200}%
\put(-\value{x},0){\oval(\value{y},#200)[l]}%
\put(-\value{x},0){\makebox(0,0){\begin{picture}(0,#200)%
\put(0,0){\vector(1,0){\value{z}}}%
\put(0,#200){\line(1,0){\value{z}}}\end{picture}}}%
\truex{2400}%
\put(-\value{x},0){\makebox(0,0)[r]{${#1}$}}%
\end{picture}}%

\newcommand{\Wncurve}[2]%
{\begin{picture}(0,0)%
\truex{1300}\truey{2000}\truez{200}%
\put(-\value{x},0){\oval(\value{y},#200)[l]}%
\put(-\value{x},0){\makebox(0,0){\begin{picture}(0,#200)%
\put(0,0){\line(1,0){\value{z}}}%
\put(0,#200){\vector(1,0){\value{z}}}\end{picture}}}%
\truex{2400}%
\put(-\value{x},0){\makebox(0,0)[r]{${#1}$}}%
\end{picture}}%

\newcount\SCALE%

\newcount\NUMBER%

\newcount\LINE%

\newcount\COLUMN%

\newcount\WIDTH%

\newcount\SOURCE%

\newcount\ARROW%

\newcount\TARGET%

\newcount\ARROWLENGTH%

\newcount\NUMBEROFDOTS%

\newcounter{x}%

\newcounter{y}%

\newcounter{z}%

\newcounter{horizontal}%

\newcounter{vertical}%

\newskip\itemlength%

\newskip\firstitem%

\newskip\seconditem%

\newcommand{\printarrow}{}%

\newcommand{\truex}[1]{%
\NUMBER=#1%
\multiply\NUMBER by 100%
\divide\NUMBER by \SCALE%
\setcounter{x}{\NUMBER}}%

\newcommand{\truey}[1]{%
\NUMBER=#1%
\multiply\NUMBER by 100%
\divide\NUMBER by \SCALE%
\setcounter{y}{\NUMBER}}%

\newcommand{\truez}[1]{%
\NUMBER=#1%
\multiply\NUMBER by 100%
\divide\NUMBER by \SCALE%
\setcounter{z}{\NUMBER}}%

\newcommand{\changecounters}[1]{%
\SOURCE=\ARROW%
\ARROW=\TARGET%
\settowidth{\itemlength}{#1}%
\ifdim \itemlength > 2800\unitlength%
\addtolength{\itemlength}{-2800\unitlength}%
\TARGET=\itemlength%
\divide\TARGET by 1310%
\multiply\TARGET by 100%
\divide\TARGET by \SCALE%
\else%
\TARGET=0%
\fi%
\ARROWLENGTH=5000%
\advance\ARROWLENGTH by -\SOURCE%
\advance\ARROWLENGTH by -\TARGET%
\advance\SOURCE by -\TARGET}%

\newcommand{\initialize}[1]{%
\LINE=0%
\COLUMN=0%
\WIDTH=0%
\ARROW=0%
\TARGET=0%
\changecounters{#1}%
\renewcommand{\printarrow}{#1}%
\begin{center}%
\vspace{10pt}%
\begin{picture}(0,0)}%

\newcommand{\DIAGV}[2]{%
\SCALE=#1%
\setlength{\unitlength}{655sp}%
\multiply\unitlength by \SCALE%
\divide\unitlength by 100%
\initialize{\mbox{$#2$}}}%

\newcommand{\n}[1]{%
\changecounters{\mbox{$#1$}}%
\put(\COLUMN,\LINE){\makebox(0,0){\printarrow}}%
\thinlines%
\renewcommand{\printarrow}{\mbox{$#1$}}%
\advance\COLUMN by 4000}%

\newcommand{\nn}[1]{%
\put(\COLUMN,\LINE){\makebox(0,0){\printarrow}}%
\thinlines%
\ifnum \WIDTH < \COLUMN%
\WIDTH=\COLUMN%
\else%
\fi%
\advance\LINE by -4000%
\COLUMN=0%
\ARROW=0%
\TARGET=0%
\changecounters{\mbox{$#1$}}%
\renewcommand{\printarrow}{\mbox{$#1$}}}%

\newcommand{\conclude}{%
\put(\COLUMN,\LINE){\makebox(0,0){\printarrow}}%
\thinlines%
\ifnum \WIDTH < \COLUMN%
\WIDTH=\COLUMN%
\else%
\fi%
\setcounter{horizontal}{\WIDTH}%
\setcounter{vertical}{-\LINE}%
\end{picture}}%

\newcommand{\diag}{%
\conclude%
\raisebox{0pt}[0pt][\value{vertical}\unitlength]{}%
\hspace*{\value{horizontal}\unitlength}%
\vspace{10pt}%
\end{center}%
\setlength{\unitlength}{1pt}}%

\newcommand{\diagv}[3]{%
\conclude%
\NUMBER=#1%
\rule{0pt}{\NUMBER pt}%
\hspace*{-#2pt}%
\raisebox{0pt}[0pt][\value{vertical}\unitlength]{}%
\hspace*{\value{horizontal}\unitlength}
\NUMBER=#3%
\advance\NUMBER by 10%
\vspace*{\NUMBER pt}%
\end{center}%
\setlength{\unitlength}{1pt}}%

\newcommand{\N}[1]%
{\raisebox{0pt}[7pt][0pt]{$#1$}}%

\newcommand{\crosslength}[2]{%
\settowidth{\firstitem}{#1}%
\settowidth{\seconditem}{#2}%
\ifdim\firstitem < \seconditem%
\itemlength=\seconditem%
\else%
\itemlength=\firstitem%
\fi%
\divide\itemlength by 2%
\hspace{\itemlength}}%

\newtheorem{definition}{Definition}[section]
\newtheorem{remark}{Remark}[section]
\newtheorem{remarks}[remark]{Remarks}
\newtheorem{notation}{Notation}{\bf}{\it}
\newtheorem{proposition}{Proposition}[section]
\newtheorem{theorem}{Theorem}[section]
\newtheorem{lemma}{Lemma}[section]
\newtheorem{corollary}{Corollary}[section]

\def\Hom{\mbox{\rm Hom}}
\def\Ext{\mbox{\rm Ext}}
\def\Tor{\mbox{\rm Tor}}

\def\D{{\rm \bf D}}
\def\J{{\rm \bf J}}
\def\Mod{{\sf -Mod}}
\def\MMod{ {\mathfrak M}(A)}
\def\Firm{{\mathfrak F}(A)}
\def\dMMod{ {\mathfrak M}(A^{op})}

\def\N{\mathbb{N}}
\def\Z{\mathbb{Z}}
\def\Lt{\mathcal{L}}
\def\S{\mathcal{S}}
\def\F{\mathcal{F}}
\def\C{\mathcal{C}}
\def\A{\mathcal{A}}
\def\FP{{\mathcal{F}}^{\perp}}

\def\li{{\displaystyle \lim_{\rightarrow}}\ }

\def\Im{\mbox{\rm Im }}
\def\Ker{\mbox{\rm ker}}
\def\Coker{\mbox{\rm coker}}
\def\dgC{{\it dg}\widetilde{{\mathcal C}}}
\def\barF{{\widetilde{{\mathcal F}}}}

\def\dgF{{\it dg}\widetilde{{\mathcal F}}}

\def\barC{{\widetilde{{\mathcal C}}}}

\newcommand{\Ch}{\mathbb{C}}

\begin{document}
\title[Flat Model Structures for Nonunital Algebras]
{Flat Model Structures for Nonunital Algebras and Higher K-Theory}

\author{Sergio Estrada and Pedro A. Guil Asensio}
\address{
Departamento de Matem\'atica Aplicada, Universidad de Murcia, Campus
del Espinardo, Espinardo (Murcia) 30100, Spain}
\email{sestrada@um.es}%
\address{
Departamento de Matem\'aticas, Universidad de Murcia, Campus del
Espinardo, Espinardo (Murcia) 30100, Spain}
\email{paguil@um.es}%

\thanks{Work partially supported by the DGI and by the Fundaci\'on S\'eneca. Part of the sources of both institutions come from the FEDER funds of the European Union.}
\subjclass[2000]{Primary: 18G55. Secondary: 16E40, 18E30, 19D55}%
\keywords{Nonunital algebras, cotorsion pairs, h-unitary modules, model structures, Waldhausen categories.}
\maketitle


\begin{abstract}
We prove the existence of a Quillen Flat Model Structure in the
category of unbounded complexes of h-unitary modules over a
nonunital ring (or a $k$-algebra, with $k$ a field). This model
structure provides a natural framework where a Morita-invariant
homological algebra for these nonunital rings can be developed.
And it is compatible with the usual tensor product of complexes.
The Waldhausen category associated to its cofibrations
allows to develop a Morita invariant excisive higher $K$-theory for nonunital
algebras.
\end{abstract}



\section{Introduction.}

Let $A$ be a nonunital algebra (or ring). A classical question in Homotopy Theory
is to find a 'good definition' of $K$-theory and cyclic type
homology for this type of rings and algebras. Namely, it is always possible to
embed a nonunital ring $A$ as a two-sided ideal of a unital ring
$R$ (for instance, by choosing $\tilde{A}=\mathbb{Z}\ltimes A$ to be the ring
obtained by adjoining an identity to $A$). Thus, it is possible to
define the notions of $K$-theory and cyclic type homologies for $A$
in terms of this ring $\tilde{A}$. But this embedding of $A$ into a
unital ring $R$ is not unique. And therefore, the different choices
of $R$ give rise to different definitions of homology theories and
$K$-theory for $A$. This problem is known in the literature as the
'excision problem' in the different theories. In \cite{wodzicki}, Wodzicki proved
that if $R$ is a unitary
$k$-pure extension of $A$, then it satisfies the excision property
for Cyclic, Bar or Hochschild homology if and only if $A$ is an
$H$-unital $k$-algebra, in the sense that its Bar homology
$HB_*(A,V)=0$ for any $k$-module $V$ (see e.g. \cite[Theorem 3.1]{wodzicki}). This result extended to
algebras over commutative rings his remarkable result showing that
if a (nonunital) ring $A$ satisfies the excision property in
rational algebraic $K$-theory, then the $\mathbb Q$-algebra
$A\otimes \mathbb{Q}$ verifies the excision property in Cyclic
homology \cite{wodzicki} (see also \cite{loday,sw}). But this $H$-unitality
condition is rarely satisfied in practice and thus, many authors
have extended these ideas and techniques to more general settings.
For instance, Cuntz and Quillen have proven in \cite{CQ1,CQ2,CQ3}
that arbitrary extensions of $k$-algebras ($k\supset \mathbb{Q}$ a
field) with $k$-linear section satisfy excision in periodic cyclic
cohomology. And Weibel in \cite{weibel2} constructed a homotopy invariant
algebraic K-theory satisfying excision and cohomologic descent.
The obstruction for the classical K-theory excision has also been studied in
\cite{cortinas} (see also \cite{cortinas2}).


In the present paper, we addopt the definition
of $h$-unitary modules given by Suslin and
Wodzicki (cf. \cite[7.3(ii)]{sw}). And we define the category $\MMod$ of
h-unitary left $A$-modules as the full subcategory of all left $A$-modules
consisting on those modules $M$ satisfying
that $\Tor_n^{\tilde{A}}(\mathbb{Z},M)=0$ (or $\Tor_n^{\tilde{A}}(k,M)=0$, if $A$ is an algebra over a ground field $k$). At first sight, this definition of $\MMod$ depends on the choosing embedding of $A$ as a two-sided ideal of the unital ring $R$. And therefore we will refer to it as the h-unitary module category associated to the pair $(R,A)$. But we will show in the last section of the paper that this h-unitary module category is independent of $R$.

We prove in Theorem
\ref{aproxuni} that any unitary left
$A$-module $M$ (in particular, the regular module $A$) has a (unique up to isomorphism) cover by an
h-unitary module (equivalently, a minimal right approximation
in the sense of Auslander \cite{AS}) by an h-unitary module. And then, we show that there is a very
satisfactory version of relative homological algebra in the category
of unbounded complexes of h-unitary modules. Namely, we prove in Section \ref{modstr} that the class of all flat h-unitary $A$-modules imposes a cofibrantly generated Model Structure (see \cite{Quillen}) in the category $\mathbb{C}(\MMod)$ of all complexes in $\MMod$.
Let us note that, if we consider in $\mathbb{C}(\tilde{A})$ the cofibrantly generated Flat Model Structure constructed in
\cite{Gill}, then the bounded flat complexes which are generating cofibrations form a small Waldhausen category $\mathcal S$ in the sense of \cite{weibel}. And the subset $\mathcal T\subset\mathcal S$ of generating cofibrations which are flat $h$-unitary is a Waldhausen subcategory of it, since any cofibration in our model structure on $\mathbb{C}(\MMod)$ is a cofibration in $\mathbb{C}(\tilde{A})$. Therefore we can use the general construction of $K$-theory on a Waldhausen category (see \cite[Chapter IV]{weibel}) to define a $K$-theory on $A$. This $K$-theory is going to be excisive in the sense of \cite{wodzicki} by our results in the last section. In particular, we deduce that the group $K_0\mathcal T$ is independent of the choice of $\mathcal S$ (see \cite[Chapter II, Theorem 9.2.2]{weibel}). The explicit construction of these $K$-groups in terms of our model structure will be developed in a forthcoming paper. We would like to
stress that our construction does not restrict just to $H$-unital algebras, but it applies to any nonunital algebra (or ring) without
the h-unitality condition.

Unfortunately, the category of h-unitary modules is far from being a monoidal category
and therefore, we cannot expect to construct a monoidal model structure on $\mathbb{C}(\MMod)$.
We solve this problem in Section 4, where we consider the wider subcategory $\Firm$ of
$\tilde{A}$-Mod consisting of all firm modules in the sense of \cite{Quillen2}.
The usual tensor product of $A$-modules is an endofunctor in $\Firm$. And we show in
Theorem \ref{escom} that this tensor product in $\Firm$ induces a unitless monoidal
structure in $\mathbb{C}(\Firm)$ which is compatible with the model structure we have
induced in $\mathbb{C}(\MMod)$.

At this point, the question of when the additive category of
h-unitary modules is abelian naturally arises. We show in Section 2 that this category is always an
accessible category in the sense of \cite{Adamek,mp}.
Actually, the constructions developed in this section will be
critical for proving our main results.

Let us finally remark that, although the category of h-unitary modules over nonunital rings extends the category of
(left) $R$-modules over a unital ring $R$ (actually, if $R$
is unital the three categories do coincide), there are important problems
for developing a notion of homological algebra in these categories. The main reason is that neither of these categories contains in general enough
projectives. In order to solve this problem, we prove that both categories have enough flat modules (note that the notion of firm and h-unitary modules do coincide in the case of flat modules). That is, that any
h-unitary module (resp. firm module) is the homomorphic image of an h-unitary flat module (resp. a firm flat module). Thus
we can define flat resolutions of h-unitary or firm modules.
However, it is not possible to define a good homotopy relation
among these flat resolutions in the category of unbounded chain
complexes. And this means that these resolutions do
not allow to uniquely define torsion functors in the corresponding
homotopy categories. We solve this problem
by introducing a good homotopy theory in the chain complexes category
of h-unitary modules.
Namely, we prove that it is possible to impose an h-unitary flat
Quillen Model Structure (cf. \cite{Hov,Quillen}) in the category of
h-unitary modules. Our proof of the existence of this Quillen Model
Structure uses Hovey's criteria (cf. \cite{hovey2}) relating
Cotorsion Pairs and Model Category Structures. The flatness condition in this model structure ensures its compatibility respect to the monoidal structure in the larger category $\Ch(\Firm)$ induced by the endofunctor given by the usual tensor product of complexes. Moreover
we deduce that h-unitary flat modules are preserved under
equivalences given by a Morita context. As an application to this construction, if we take the nonunital algebra to be a $C^*$-algebra, our model category structure provides a solution to a question posed by Hovey in \cite[Problem 8.4]{Hov}.

Along this paper, all rings will be associative and non necessarily
unitary. Although we state our results for nonunital rings $A$, we
can assume (as it often occurs in practice) that our nonunital rings are
algebras over some field $k$ (so $\tilde{A}$ is the $k$-algebra
$k\ltimes A$ obtained by adjoining an identity to $A$). We refer to
\cite{GoTr,Hov,loday, weibel} for any undefined notion on cotorsion pairs,
Model Structures or K-theory used in the text.

\section{Main Lemmas.}

Let $\mathcal A$ be a Grothendieck category and $\Lt$, a class of
objects of $\mathcal A$ closed under isomorphisms. A well-ordered direct system
 of objects of $\mathcal A$ $( X_{ \alpha } | \alpha \leq \lambda ) $
is said to be continuous if $ X_0  = 0, X_{\alpha}\subseteq
X_{\alpha+1} $ and, for each limit ordinal $ \beta \leq \lambda $,
we have that $ \displaystyle X_{\beta}=\cup_{\alpha<\beta} X_{
\alpha } $. An object $X$ of $\mathcal A$ is called $\Lt$-filtered
provided that $ X = X_{ \lambda } $ for some well-ordered direct system $ ( X_{ \alpha } | \alpha \leq \lambda ) $ such that
$ X_{ \alpha+1 }/ X_{ \alpha }$ is in $ {\mathcal L} $ whenever
$\alpha + 1 \leq \lambda $.

\begin{notation}
Let ${\mathcal D}$ be a class of objects of a Grothendieck category
${\mathcal{A}}$. We will denote by  ${\mathcal D}^{\perp}$ the class
of all objects $Y$ of ${\mathcal{A}}$ such that $\Ext^1(D,Y)=0$ for
every $D\in {\mathcal D}$. Similarly, $^{\perp}{\mathcal D}$ will
denote the class of those objects $Z\in\mathcal{A}$ such that
$\Ext^1(Z,D)=0$ for every $D\in {\mathcal D}$.
\end{notation}

We recall that a pair (${\mathcal F}, {\mathcal C}$) of classes of
objects of ${{\mathcal{A}}}$ is called a cotorsion pair if
${\mathcal F}^{\perp}={\mathcal C}$ and $^{\perp}{ \mathcal
C}={\mathcal F}$ (see e.g. \cite{GoTr}). The cotorsion pair is said
to have enough injectives (resp. enough projectives) if, for every
$Y$ in ${{\mathcal{A}}}$, there exists an exact sequence
 $0\rightarrow Y\rightarrow C\rightarrow
F\rightarrow 0$ (resp. for every $Z\in\mathcal{A}$ there exists an
exact sequence $0\rightarrow C'\rightarrow F'\rightarrow
Z\rightarrow 0$) where $F,F'\in {\mathcal F}$ and $C,C'\in {\mathcal
C}$. The cotorsion pair $({\mathcal F},{\mathcal C})$ is complete if
it has enough injectives and projectives. We will say that
$({\mathcal F}, {\mathcal C}) $ is functorially complete when these
sequences can be chosen in a functorial manner (depending on $Y$ and
$Z$) (see \cite[Definition 2.3]{hovey2}). Finally, a cotorsion pair
$(\F,\FP)$ is said to be cogenerated by a set $\mathcal
S\subseteq\mathcal{F}$ if $\mathcal S^{\perp}=\FP$. From results in
\cite{Estr2} we get the following Theorem.

\begin{theorem}\label{bSalce}
Let $\F$ be a class of objects of a Grothendieck category $\A$ which
is closed under direct sums, extensions and well ordered direct
limits. Suppose that the pair $(\F,\FP)$ is cogenerated by a set. If
$\F$ contains a generator of $\A$, then the pair $(\F,\FP)$ is a
complete cotorsion pair.
\end{theorem}

\begin{proof} It follows from \cite[Theorem 2.5]{Estr2} that the pair
$(\F,\FP)$ has enough injectives. To show that the pair also has
enough projectives, we will adapt the arguments of
\cite[Lemmas 2.2 and 2.3]{Salce}. Let $M$ be any object of $\A$.
There exists a short exact sequence $0\to K\to G\to M\to 0$ with
$G\in \F$. Since $(\F,\FP)$ has enough injectives there also exists
a short exact sequence $0\to K\to C\to F\to 0$, with $F\in \F$ and
$C\in \FP$. Let us construct the pushout diagram
 \DIAGV{50} {} \n{} \n{0}\n{}\n{0}\nn
{}\n{}\n{\sar}\n{}\n{\sar}\nn{0}
\n{\ear}\n{K}\n{\ear}\n{G}\n{\ear}\n{M}\n{\ear}\n{0}\nn
{}\n{}\n{\sar}\n{}\n{\sar}\n{}\n{\seql}\nn
{0}\n{\ear}\n{C}\n{\ear}\n{Q}\n{\ear}\n{M}\n{\ear}\n{0}\nn
{}\n{}\n{\sar}\n{}\n{\sar}\nn {}\n{}\n{F}\n{\eeql}\n{F}\nn
{}\n{}\n{\sar}\n{}\n{\sar}\nn {}\n{}\n{0}\n{}\n{0}\diag and let us
note that $Q\in \F$, since $\F$ is closed under extensions. Thus,
the short exact sequence $0\to C\to Q\to M\to 0$ shows that
$(\F,\FP)$ has enough projectives, i.e. it is a complete pair.
Finally if $M\in \!^{\perp}\F$ then the previous short exact sequence
splits. In particular, $M\in\F$ and hence $(\F,\FP)$ is a cotorsion pair.
\end{proof}

Let us fix a nonunital ring $A$. And let us denote by
$\tilde{A}=\Z\ltimes A$ the unital ring obtained from $A$ by
adjoining an identity. A left $A$-module $M$ is said to be h-unitary
if the evaluation map $\mu:A\otimes_{\tilde{A}} M\to M$ given by
$\mu(a\otimes m)=am$ is an isomorphism and
$\Tor_i^{\tilde{A}}(A,M)=0, \forall i\geq 1$ or, equivalently, if
$\Tor_n^{\tilde{A}}(\Z,M)=0$, $\forall n\geq 0$. Let us note that
this definition generalizes the corresponding one given in \cite[7.2
(ii)]{sw} for modules over nonunital $k$-algebras. Results in
Section \ref{moritinv} show that the category of h-unitary modules
remains invariant if we replace $\tilde{A}$ by any other unitary
extension of the ring $A$ (i.e., by any other unitary ring $R$
containing $A$ as a two-sided ideal). Therefore, from now on, {\em
we will denote by $R$ any unitary extension ring of $A$} and we will
implicitly assume, if necessary, that $R=\tilde{A}$. Other
concept which is independent of the chosen unitary extension ring of $A$ is that
of {\it firm} module. A left $A$-module $M$ is said to be firm if
$\mu$ is an isomorphism or, equivalently, if $\Tor_i^{\tilde{A}}(\Z,M)=0$ for
$i=0,1$. We denote by $\MMod$ and $\dMMod$ (resp. $\Firm$ and
$\mathfrak{F}(A^{op})$) the categories of left and right
h-unitary $A$-modules (resp. left and right firm modules). We
summarize next the main properties of these subcategories of $R\Mod$.
We recall that a morphism in a category is called a monocokernel
(resp., epikernel) if it is the cokernel (resp. kernel) of a monomorphism (resp., epimorphism).

\begin{lemma}\label{genpl}
$\MMod\subseteq \Firm$ are both additive subcategories of $R\Mod$
having a flat generator. The category $\Firm$ is cocomplete. The
category $\MMod$ is closed under monocokernels, epikernels and
extensions in $R\Mod$.
\end{lemma}
\begin{proof}
Let us check the first claim. We prove the following slightly more
general statement: let $M$ be a unitary left $A$-module (that is
$AM=M$). Then $M$ is the homomorphic image of a unitary flat
$A$-module. To see this claim, let us fix an epimorphism $g:P\to M$
from a projective $R$-module $P$ onto $M$. The restriction $g|_{AP}$
is also surjective (since $AM=M$) and hence, there exists an
endomorphism $f:P\to P$ such that $g\circ f=g$ and $f(P)\subseteq
AP$. Let $F=\li_{n< w_0} P_n$, where $P_n=P$ and $P_n\to P_{n+1}$ is
$f$, $\forall n\geq 0$. $F$ is a flat $R$-module by construction and
$AF=F$ since $f(P_n)\subseteq AP_{n+1}$. In particular, $F$ is an
h-unitary module. As $g:P_n\to M$ is surjective $\forall n\geq 0$,
the induced map $F\to M$ is also surjective.

The category $\Firm $ is clearly cocomplete because it is closed under
coproducts (since coproducts commute with tensor products) and it is easy to check that
for every morphism $f:M\to N$, $\Coker(f)$ computed in $R\Mod$ is actually in
$\Firm$.

Finally, let $0\to M_1\to M_2\to M_3\to 0$ be any short exact sequence in
$R\Mod$ such that $M_i, M_j\in \MMod$, for $i, j\in \{1,2,3\},\ i\neq
j$. Then by taking the long exact sequence with respect to
$\Tor^R_*(\Z,?)$ we get that $M_k\in \MMod$, with $k\in \{1,2,3\},\
k\neq i,j$.
\end{proof}

\begin{remarks}

\begin{enumerate}
\item [(1)] Let us note that the category $\Firm$ is also closed under extensions and epikernels in $R\Mod$.

\item[(2)]We do not know whether the category $\MMod$ is cocomplete in general.
\end{enumerate}
\end{remarks}

Our next lemma extends to sets of arbitrary cardinality a well-known
result that belongs to the folklore of the Theory of Purity of
modules. It has been used, for instance, by Bican, El Bashir and
Enochs for proving the existence of flat covers of modules (see
\cite{BiBash}). We denote, as usual, the cardinality of a set $N$ by
$\left\vert N\right\vert$.

\begin{lemma}
Let $R$ be a unitary ring and $M$, a unitary left $R$-module. Every
subset $N$ of $M $ is contained in a pure submodule of $M$ of
cardinality at most $\aleph =max\left\{ \left\vert R\right\vert
,\left\vert N\right\vert ,\aleph _{0}\right\} $.
\end{lemma}

\begin{proof}

Let $N_{0}$ be the submodule of $M$ generated by $N$. Note that
$N=\sum_{y\in N}Ry$ and therefore, $\vert N_0\vert \leq \sum_{y\in
N} \vert Ry\vert \leq \vert N\vert\cdot\vert R\vert=\aleph$.

Consider now the short exact sequence of modules
\begin{equation*}
0\rightarrow N_{0}\overset{u}{\longrightarrow }M\overset{p}{\longrightarrow }%
M/N_{0}\rightarrow 0
\end{equation*}
and let $\Omega _{0}$ be the subset of $\cup _{n,m\in
\mathbb{N}}\left[ \Hom_{R}(R^{n},N_{0})\times
\Hom_{R}(R^{n},R^{m})\right] $ consisting on those
pairs of morphisms $\left( h,v\right) $ such that there exists a $%
g_{(h,v)}:R^{m}\rightarrow M$ with $g_{(h,v)}\circ v=u\circ h$. Note
that
\begin{equation*}
\begin{array}{c}
\left\vert \Omega _{0}\right\vert \leq \left\vert \cup _{n,m\in \mathbb{N}}%
\left[ \Hom_{R}(R^{n},N_{0})\times \Hom_{R}(R^{n},R^{m})\right]
\right\vert =\\ \sum\nolimits_{n,m\in \mathbb{N}}\left\vert
\Hom_{R}(R^{n},N_{0})\times \Hom_{R}(R^{n},R^{m})\right\vert =
\\
\sum\nolimits_{n,m\in \mathbb{N}}\left\vert N_{0}^{n}\times
R^{n\times m}\right\vert \leq max\left\{ \left\vert R\right\vert
,\left\vert
N\right\vert ,\aleph _{0}\right\}.%
\end{array}%
\end{equation*}%
And let $N_{1}$ be the submodule of $M$ generated by
$$N_{0}\cup \left\{ \Im(g_{(h,v)})\mid \,(h,v)\in \Omega _{0}\right\} .$$
Then $\left\vert \Im(h)\right\vert \leq max\left\{ \left\vert
R\right\vert ,\left\vert N\right\vert ,\aleph _{0}\right\} $ for any
$(h,v)\in \Omega _{0} $ and therefore, $N_{1}$ has at most $\aleph
\times \aleph =\aleph $ generators. In particular, this means that again
$\left\vert N_{1}\right\vert \leq \aleph $, since $\aleph \geq
\left\vert R\right\vert $.

Let us now replace $N_{0}$ by $N_{1}$ and let us construct $N_{2}$ in a
similar way. Following this method, we will get, by induction on
$\mathbb{N}$, an infinite ascending chain $\left\{ N_{k}\right\}
_{k\in \mathbb{N}}$ of submodules of $M$ of cardinality bounded by
$\aleph $. Let $L=\cup _{k\in \mathbb{N}}N_{k}$. Clearly $\left\vert
L\right\vert$ is also bounded by $\aleph $.

We claim that $L$ is pure in $M$. To check it, let $\phi
:F\rightarrow M/L$ be a homomorphism from a finitely presented
module to $M/L$ and let us show that there exists a morphism
$t:F\rightarrow M$ such that $\phi=p\circ t$, where $p:M\rightarrow
M/L$ is the canonical projection. $F$ is always the cokernel of a
homomorphism $v:R^{n}\rightarrow R^{m}$. By projectivity, there
exist homomorphisms $g:R^{m}\rightarrow M$ and $h:R^{n}\rightarrow
L$
such that the following diagram commutes%

\medskip\par\noindent
\begin{equation*}
\begin{array}{ccccccccc}
&  & R^{n} & \overset{v}{\longrightarrow } & R^{m} & \overset{q}{%
\longrightarrow } & F & \rightarrow & 0 \\
&  & \downarrow ^{h} &  & \downarrow ^{g} &  & \downarrow ^{\phi } &  &  \\
0 & \rightarrow & L & \overset{u}{\longrightarrow } & M & \overset{p}{%
\longrightarrow } & M/L & \rightarrow & 0%
\end{array}%
\end{equation*}

\medskip\par\noindent
As $\Im(h)$ is a finitely generated submodule of $L=\cup _{k\in \mathbb{N}%
}N_{k}$, there exists a $k\in \mathbb{N}$ such that $\Im(h)\subset
N_{k}$. Therefore, the element $(h,v)\in \Omega _{k+1}$ and, by
construction, this means that there exists a homomorphism $g^{\prime
}_{(h,v)}:R^{m}\rightarrow N_{k+1}\subset L$ such that $g^{\prime
}_{(h,v)}\circ v=h$. But then $\left( g-g^{\prime }_{(h,v)}\right)
\circ v=0$ and thus, there exists a $t:F\rightarrow M$ such that
$t\circ q=g-g^{\prime }_{(h,v)}$. In particular, $p\circ t\circ
q=p\circ (g-g^{\prime }_{(h,v)})=\phi \circ q-p\circ u\circ
g^{\prime }_{(h,v)}=\phi \circ q$. Therefore, $\phi =p\circ t$,
since $q$ is an epimorphism.  \end{proof}

Let us now consider the category $\MMod$ of left h-unitary modules.
Let us recall that, by Lemma \ref{genpl}, $\MMod$ has an h-unitary
flat generator $G$.

\begin{lemma}\label{tip}
Let $G$ be an h-unitary flat generator of $\MMod$ and let us fix a
cardinal number $\aleph \geq max\left\{ \left\vert G\right\vert,
\aleph _{0}\right\}$. Let $M$ be an h-unitary  $A$-module and $N$, a
subset of $M$ of cardinality bounded by $\aleph$. Then $N$ embeds
in a pure h-unitary submodule $L$ of $M$ such that $\left\vert L\right\vert \leq \aleph$.
\end{lemma}

\begin{proof}
Let us denote by $\rho :A\otimes _{R}M\rightarrow M$ the structural
homomorphism. Note that $\rho$ is an isomorphism since $M$ is
h-unitary. Let us construct $L$ by induction on $n\in \mathbb{N}$.
By the above Lemma, there exists a pure $R$-submodule $N_{0}\subset
M$ containing $N$ such that $\left\vert N_{0}\right\vert \leq \aleph
$. Therefore, we have the
following diagram%
\begin{equation*}
\begin{array}{ccc}
A\otimes _{R}N_{0} &  & N_{0} \\
\downarrow &  & \downarrow \\
A\otimes _{R}M & \overset{\rho }{\longrightarrow } & M
\end{array}%
\end{equation*}%
in which the vertical arrows are monomorphisms, since $N_0$ is pure
in $M$. Let us choose, for any $x\in N_{0}$, elements
$m_{1}^{x},\ldots ,m_{k_{x}}^{x}\in M$ and elements
$a_{1}^{x},\ldots ,a_{k_{x}}^{x}\in A$ such that $\rho
\left(\sum\nolimits_{i=1}^{k_{x}}a_{i}^{x}\otimes
m_{i}^{x}\right)=x$. And let us call $K_1=N_{0}\bigcup \left( \cup
_{x\in N_{0}}\left\{ m_{1}^{x},\ldots ,m_{k_{x}}^{x}\right\}
\right)$. Note that $\left\vert K_1\right\vert$ is also bounded by
$\aleph$.

On the other hand, as $G$ is a generator of $\MMod$, there exists a homomorphism $f_y:G\rightarrow M$ and an
element $z_y\in G$ such that $y=f_y(z_y)$, for any $y\in
K_1$. Let $K'_1=\sum_{y\in
K_1}f_y(G)$. Then $K'_1$ is a submodule of $M$ containing $K_1$ and thus,
also containing $N_0$. Moreover, it is straightforward to check that
$\left\vert K'_1\right\vert$ is also bounded by $\aleph$.

Let us now replace $N$ by $K'_{1}$ and continue the construction.
Set finally $L=\cup _{n\in \mathbb{N}}N_{n}=\cup _{n\in
\mathbb{N}}K'_{n}$. Clearly $L$ is a pure submodule of $M$, since it
is the union of the chain $\left\{ N_{n}\right\} _{\mathbb{N}}$ of
pure submodules of $M$. From this it is easy to deduce that
$N\subseteq M$ is also h-unitary. And $\left\vert L\right\vert \leq
\aleph _{0}\times \aleph =\aleph $.

\end{proof}

We recall that a category $\mathcal C$ is called $\aleph$-accessible
(where $\aleph$ is an infinite  regular cardinal number) if it has $\aleph$-direct
limits and there exists a set $\mathcal{C}_0$ of $\aleph$-presentable objects of
$\mathcal C$ such that
any other object of $\mathcal C$ is (isomorphic to) an
$\aleph$-direct limit of morphisms among objects in $\mathcal{C}_0$.
Where an $\aleph$-direct limit of morphisms is the direct limit of a
directed set of morphisms $\{f_{ij}\}_I$ satisfying that for any
subset $I_0\subseteq I$ of cardinality strictly smaller than
$\aleph$, there exists an $i_0\in I$ such that $i\leq i_0$ for any
$i\in I_0$. And an object $C\in\mathcal{C}$ is called $\aleph$-presentable
if the functor $\rm{Hom}(C,-):\mathcal{C}\rightarrow Set$ commutes
with $\aleph$-direct limits (see e.g. \cite{Adamek,mp}). A category is simply called accessible if it is $\aleph$-accessible for some infinite regular cardinal
$\aleph$.

\begin{corollary}\label{estec}
The category of h-unitary  left modules over a nonunital ring is an
accessible additive category.
\end{corollary}

\begin{proof}
The category is clearly additive and it has $\lambda$-direct limits,
for any infinite regular cardinal $\lambda$, since the tensor product is a
right exact additive functor commuting with direct limits. On the
other hand, the above result shows that, for any cardinal $\aleph \geq
max\left\{\left\vert G\right\vert,\aleph _{0}\right\}$, any object in the category
is the ($\aleph^+$)-direct union of its h-unitary submodules of
cardinality bounded by $\aleph$. Moreover, it is easy to check that, if we choose $\aleph\geq \vert A\vert$, the above submodules are $\aleph^+$-presentable. Finally, let
us note that $\aleph^+$ is regular, since it is a successor
cardinal.
\end{proof}

\begin{remark}

\item Let us fix a generator $G$ of $\MMod$ and let $\aleph \geq max\left\{
\left\vert G\right\vert,\aleph _{0}\right\}$. Let us call
$G_{\aleph}$ the direct sum of all isomorphism classes of h-unitary
modules of cardinality bounded by $\aleph$. The proof of the above
corollary shows that $G_{\aleph}$ is also a generator of $\MMod$.
Moreover, $G_{\aleph}$ has the following interesting property: any
h-unitary module is an ($\aleph^+$-)directed union of h-unitary
pure submodules isomorphic to direct summands of $G_{\aleph}$.
\end{remark}

\section{A Quillen Model Structure on $\Ch(\MMod)$.}\label{modstr}

The main purpose
of this section will be to impose a Quillen Model Structure (see
\cite{Hov,Quillen} for its definition) in the category $\Ch(\MMod)$
of unbounded chain complexes of h-unitary left $A$-modules in terms
of h-unitary flat modules. Our proof is based on Hovey's
Theorem \cite{hovey2} which relates Complete Cotorsion pairs and
Model Structures in the corresponding category of unbounded
complexes.

The additive category $\MMod$ may not be abelian in general, but it can be obviously
embedded as a full subcategory of the abelian category $R\Mod$. As $\MMod$
is closed under extensions in $R\Mod$ (Lemma
\ref{genpl}), we may consider the proper class  $\mathcal{P}$ in $\MMod$
consisting of all short
exact sequences in $R\Mod$, $0\to M\to N\to P\to 0$ with $M,N,P\in
\MMod$. And then define as in (\cite[Chapter XII, Theorem 4.4]{maclane}), the relative extension groups
$\Ext^n_{\MMod}(M,N)$,  for all $M,N\in \MMod$ and $n\geq 0$, as well as natural transformations of bifunctors $\Ext_{\MMod}^n
(?,?)\to \Ext_R^n(?,?)$, for all $n\geq 0$. Note that these natural transformations are monomorphisms
for $n=1$, and isomorphisms for $n=0$.

We recall from the above section that an h-unitary module $A$ is flat
iff it is flat in $\tilde{A}\Mod$. Let us denote by $\F$ the class
of all h-unitary flat $A$-modules. We will show in our last section
(Theorem \ref{esMi}) that this class does not depend on the chosen embedding
of $A$ into a unital ring $R$. Let us also note that a flat
$R$-module $F$ belongs to $\F$ if and only if $A\cdot F=F$

Let us denote by $\C$ the class of all h-unitary modules which are
Ext-orthogonal to $\F$, that is,
$$\C=\{C\in \MMod:\ \Ext^1_{\MMod}(F,C)=0, \forall F\in \F\}.$$
Modules in $\C$ will be called h-unitary cotorsion modules.

The problem of the previous approach is that these relative Ext groups are computed in
terms of extensions in $R\Mod$.
But we would like to be able to compute them intrinsically by resolutions and
coresolutions in $\MMod$. On the other hand, it is implicit in \cite{wodzicki} that
h-unitary flat modules are the right choice for computing all kinds of
Morita invariant homology functors. So our approach will be to
use the class $\mathcal{P}$ as the proper class of exact sequences
in $\MMod$, but considering h-unitary flat and h-unitary cotorsion modules in
$\MMod$ as the basic bricks for constructing resolutions and coresolutions in $\MMod$
respectively.

From now on, when we refer to exact sequences in $\MMod$ we will mean
exact sequences in $\mathcal{P}$. We begin by proving that every
$R$-module can be approximated in a minimal fashion by a unique up
to isomorphism h-unitary module.

\begin{theorem}\label{aproxuni}
Let $M$ be any $R$-module. Then there exists an h-unitary
module $M_H$ and a morphism $\varphi:M_H\to M$ satisfying that:
\begin{enumerate}
\item For every morphism $\psi:N_H\to M$ with $N_H\in \MMod$ there
exists a morphism $f:N_H\to M_H$ such that $\varphi\circ f=\psi$.
\item Every morphism $g:N_H\to N_H$ such that $\varphi\circ g=\varphi$
is an automorphism.
\end{enumerate}
Such $\varphi:M_H\to M$ is necessarily unique up to isomorphisms and
it is called the the $\MMod$-cover of $M$.

\end{theorem}
\begin{proof}

The proof is similar to \cite[Theorem 4.1]{Est}. We just need to use Quillen's small object argument (see \cite[Lemma II.3.3]{Quillen}) and use the fact that $\MMod$ is closed under extensions, well ordered direct
limits and direct sums. Then we may apply Lemma \ref{tip} and the same
arguments of \cite[Theorem 4.1]{Est} show that every
$R$-module $M$ has an $\MMod$-cover.

\end{proof}

\begin{theorem}\label{firmcomp}
The pair $(\F,\C)$ is a complete cotorsion pair in $\MMod$.
\end{theorem}
\begin{proof} Let $F$ be an h-unitary flat module. By Lemma \ref{tip}, for
any $y\in F$ there exists a pure and h-unitary  submodule
$S\subseteq F$, with $y\in S$ and $|S|\leq \aleph$ (where $\aleph$ is the
infinite cardinal obtained in Lemma \ref{tip}). Hence, $S$ and $F/S$ are also
h-unitary flat modules. If $\S_{\aleph}$ denotes the set of all
isomorphism classes of h-unitary flat modules of cardinality bounded
by $\aleph$, we get that each h-unitary flat module $F$ is
$\S_{\aleph}$-filtered. It follows from \cite[Lemma 1]{EkTrl} that
the set $\S_{\aleph}$ cogenerates the pair $(\F,\FP)$.
Therefore by Theorem \ref{bSalce} the pair $(\F,\FP)$ is
complete. Furthermore, by Lemma \ref{genpl}, $\MMod$ is closed under
extensions, monocokernels and epikernels so it follows that in fact
the pair $(\F,\C)$ in $\MMod$ is complete.

Let us finally check that $(\F,\C)$ is a cotorsion pair. Let $M$ be
an h-unitary module such that $M\in ^{\perp}\!\!\mathcal{C}$. Since
$(\F,\C)$ has enough projectives, there exists a short exact
sequence of h-unitary modules
$$0\to K\to F\to M\to 0$$ such that $F\in \F$ and $K\in \C$.
Thus, this short exact sequence splits and $M$ is a direct summand
of $F$. In particular, $M\in \F$.

Conversely assume that $N\in \FP$ and is h-unitary. Then $\Ext^1_R(F,N)=0$, for all $F\in \F$. As $\Ext^1_{\MMod}(?,?)$ is a subfunctor of $\Ext^1_R(?,?)$ we get that
$\Ext^1_{\MMod}(F,N)=0$. Thus $N\in \C$.  \end{proof}

\bigskip\par\noindent

We will denote the category of unbounded complexes of h-unitary
modules by $\Ch(\MMod)$. That is, those complexes
$$M=\cdots\to
M^{i-1}\stackrel{\delta^{i-1}}{\longrightarrow}M^i\stackrel{\delta^{i}}{\longrightarrow}M^{i+1}\to
\cdots$$ such that $0\to \ker{\delta^{i}}\to M^i\to \Im{\delta^i}\to
0$ is in $\mathcal{P}$, $\forall i\in \Z$. We shall denote by $Z_iM$
the $i$th cycle module $\Ker(\delta^i)$ and by $B_iM$, the $i$th
boundary module $\Im(\delta^{i-1})$.

Recall that the tensor and $\Hom_R$ functors on $\MMod$ can be
canonically extended to $\Ch(\MMod)$ as follows: If $M$ and $N$ are
complexes of h-unitary modules, we call $Hom(M,N)$ the complex of
abelian groups satisfying that
$$Hom(M,N)^n=\prod_{t\in \Z}\Hom_{R}(M^t,N^{n+t})$$
and such that if $f\in Hom(M,N)^n$, then
$$(\delta^nf)^m=\delta_N^{m+n}\circ f^m-(-1)^nf^{m+1}\circ
\delta_M^m.$$ The tensor product of a complex of right h-unitary
modules $M$ and a complex of left h-unitary modules $N$ is the
complex of abelian groups $M\otimes N$ with $(M\otimes
N)^n=\oplus_{t\in \Z}M^t\otimes_{R}N^{n-t}$ and
$$\delta(x\otimes y)=\delta^t_M(x)\otimes y+(-1)^t x\otimes
\delta_N^{n-t}(y),$$for $x\in M^t$ and $y\in N^{n-t}$.

$\Ch(\MMod)$ is an additive category, since so is $\MMod$. A complex
$M$ such that $Z_iM=B_iM$, $\forall i\in \Z$ will be called {\it
exact}. We shall denote by $\mathcal{E}$ the class of all exact
complexes. If $G$ is an h-unitary flat generator of $\MMod$ then the
family of complexes $\{S^n(G):n\in \Z\}$ generates the category
$\Ch(\MMod)$. Where $S^n(G)$ denotes the complex with $G$ in the $n$'th
position and 0 in all other places. If $M$ is an h-unitary module we
will denote by $D^n(M)$ the complex $\cdots\to 0\to
M\stackrel{id}{\to}M\to 0\to \cdots$ with $M$ in positions $n-1$ and $n$.

We recall from \cite{Gill} the following definitions. A complex of
h-unitary modules $C$ is said to be {\it $\F$-cotorsion} if $C$ is
exact and  $Z_nC\in \C$, for all $n\in \Z$. A complex $F$ is called
a dg-h-unitary flat complex if $F^n$ is an h-unitary flat module for
any $n\in \Z$ and, for every $\F$-cotorsion complex $C$, $Hom(F,C)$
is an exact complex in $\Z$-Mod. Dually, we define h-unitary flat
complexes as those complexes $X$ such that $X$ is exact and $Z_nX\in
\F$, $\forall n\in \Z$. And dg-cotorsion complexes, as those $Y$
such that $Y_n\in \C$ and $Hom(X,Y)$ is exact, for any h-unitary
flat complex $X$. We shall denote by $\barF$ and $\barC$ the classes
of h-unitary flat complexes and $\F$-cotorsion complexes,
respectively. And by $\dgF$ and $\dgC$, the classes of dg-h-unitary
flat and dg-cotorsion complexes of h-unitary modules. It is clear
that if $F$ is an h-unitary flat module, then $S^n(F)\in \dgF$ for
all $n\in \Z$. In particular, $\dgF$ contains the previous family of
generators of $\Ch(\MMod)$.

In order to obtain a Model Structure in $\Ch(\MMod)$, we are going
to apply Hovey's Theorem relating Cotorsion Pairs and Model Category
Structures.

\begin{theorem}\cite[Theorem 2.2]{hovey2}\label{princiHov}
Let $\Ch(\MMod)$ be the category of chain complexes of h-unitary
modules. Let $\mathcal{E}$ be the class of exact complexes in
$\MMod$. If $(\dgF,\dgC\cap\mathcal{E})$ and $(\dgF\cap \mathcal{E},
\dgC)$ are complete cotorsion pairs, then there exists a Model
Structure on $\Ch(\MMod)$. In this Model Structure, the weak
equivalences are the homology isomorphisms, the cofibrations are the
monomorphisms whose cokernels are in $\dgF$, and the fibrations are
the epimorphisms whose kernels are in $\dgC$.
\end{theorem}

\begin{remark}
The above theorem is proved in \cite[Theorem 2.2]{hovey2} under the assumption that the considered category is Abelian. We do not know if the category $\MMod$ is Abelian, but one can check that it satisfies all the conditions used in the proof of \cite[Theorem 2.2]{hovey2}. Namely, it is an additive category having a generator and it is closed under extensions, monocokernels, epicokernels and direct limits.
\end{remark}

In order to apply the above criterium, we will need to prove the
following:
\begin{enumerate}
\item The pairs $(\barF,\dgC)$ and $(\dgF,\barC)$ are cotorsion
pairs.
\item The pairs $(\barF,\dgC)$ and $(\dgF,\barC)$ are complete.
\item $\dgF\cap \mathcal{E}=\barF$ and $\dgC\cap \mathcal{E}=\barC$
\end{enumerate}

Before proving all these conditions, we need to fix the
statements of \cite[Lemma 3.8(7),(8)]{Gill}.
\begin{lemma}\label{mist}
Let $M,N$ be two
complexes of $\Ch(\MMod)$ and $C$, an object of $\MMod$.
Then there exist canonical monomorphisms of abelian groups
$$0\to\Ext^1_{\MMod}(C,Z_nN) \to \Ext^1_{\Ch(\MMod)}(S^n(C),
N)$$ and $$0\to\Ext^1_{\MMod}(M_n/B_nM,C) \to
\Ext^1_{\Ch(\MMod)}(M,S^{n}(C)).$$
\end{lemma}
\begin{proof} Let
$$0\to Z_nN\to T\to C\to
0$$ be any extension in $\Ext^1_{\MMod}(C,Z_n N)$, and let us
construct the pushout of the inclusions $Z_n N\to N_n$ and $Z_n N\to
T$ \DIAGV{60} {} \n{} \n{0}\n{}\n{0}\nn
{}\n{}\n{\sar}\n{}\n{\sar}\nn{0} \n{\ear}\n{Z_n
N}\n{\ear}\n{T}\n{\ear}\n{C}\n{\ear}\n{0}\nn
{}\n{}\n{\sar}\n{}\n{\sar}\n{}\n{\seql}\nn
{0}\n{\ear}\n{N_n}\n{\Ear{g}}\n{Q}\n{\ear}\n{C}\n{\ear}\n{0}\nn
{}\n{}\n{\sar}\n{}\n{\saR{h}}\nn
{}\n{}\n{B_{n+1}N}\n{\eeql}\n{B_{n+1}N}\nn
{}\n{}\n{\sar}\n{}\n{\sar}\nn {}\n{}\n{0}\n{}\n{0}\diag

We have a commutative diagram \DIAGV{70} {0}
\n{\ear}\n{N_{n-1}}\n{\eeql}\n{N_{n-1}}\n{\ear}\n{0}\n{\ear}\n{0}\nn
{}\n{}\n{\Sar{\delta^{n-1}}}\n{}\n{\saR{g\circ
\delta^{n-1}}}\n{}\n{\sar}\n{}\nn
{0}\n{\ear}\n{N_n}\n{\Ear{g}}\n{Q}\n{\Ear{t}}\n{C}\n{\ear}\n{0}\nn
{}\n{}\n{\Sar{\delta^n}}\n{}\n{\Sar{h}}\n{}\n{\sar}\nn
{0}\n{\ear}\n{N_{n+1}}\n{\eeql}\n{N_{n+1}}\n{\ear}\n{0}\n{\ear}\n{0}\diag
and hence, we get an extension $\xi\equiv0\to N\to H\to S^n(C)\to 0$
in $$\Ext^1_{\Ch(\MMod)}(S^n(C),N).$$ This defines a map
$$\Ext^1_{\MMod}(C,Z_n N)\to\Ext^1_{\Ch(\MMod)}(S^n(C),N).$$

Clearly this map is a morphism of abelian groups. Let us check that
it is injective. Assume that $\xi$ splits and $r:C\to Q$ is the corresponding
excision in the $n$th component of $\xi$. We get that $h\circ r=0$
by the commutativity of the diagram. Thus, ${\rm Im}(r)\subseteq T$.
Hence $0\to Z_n N\to T\to C\to 0$ splits. The proof of the second
monomorphism is dual.
\end{proof}

\begin{remark}
The monomorphisms that appear in Lemma \ref{mist} are not isomorphisms in
general, as claimed in \cite{Gill}. For instance, if we consider the category $\Ch(R)$ of
unbounded complexes of $R$-modules over a ring with identity $R$,
then $\Ext^1_R(P,-)=0$ for any projective $R$-module $P$. But the
functor $\Ext^1_{\Ch(R)}(S^n(P), -)\neq 0$, since $S^n(P)$ is a
dg-projective complex that is not projective. The same holds for our
second monomorphism if we choose an injective $R$-module.
\end{remark}

Let us now prove Condition $(1)$ in Hovey's criteria.
\begin{lemma}\label{cotpa}
The pairs $(\barF,\dgC)$ and $(\dgF,\barC)$ are cotorsion pairs in
$\Ch(\MMod)$.
\end{lemma}
\begin{proof}
The cotorsion pair $(\F,\C)$ in $\MMod$ is complete (by Theorem \ref{firmcomp}). So the result follows from \cite[Proposition
3.6]{Gill} applied to the pair $(\F,\C)$ (by realizing that the
part of the proof of \cite[Proposition 3.6]{Gill} involving
\cite[Lemma 3.8(7),(8)]{Gill} can be replaced by Lemma \ref{mist}).
\end{proof}

\bigskip
We check now Condition $(2)$.
\begin{theorem}\label{comp1}
The pair $(\dgF,\barC)$ is complete in $\Ch(\MMod)$.
\end{theorem}
\begin{proof} We already know, by the proof of Theorem \ref{firmcomp}, that the pair
$(\F,\C)$ is cogenerated by the set $\S_{\aleph}$. Let $G$ be the
h-unitary flat generator of $\MMod$. We claim that the cotorsion
pair $(\dgF, \barC)$ is cogenerated by the set
$\mathcal{L}=\{S^m(F): m\in \Z,\, F\in \S_{\aleph}\cup{\{G\}}\} $.

It is easy to check that $\mathcal{L}\subseteq \dgF$, since
$S^l(F)_l\in \F$, for any $l\in \Z$ and any $F\in
\S_{\aleph}\cup{\{G\}}$. Moreover, for every exact complex $M\in
\barC$, $Hom(S^l(F),M)$ is the complex
$$\cdots\to \Hom(F,M^l)\to
\Hom(F,M^{l+1})\to\cdots$$ which is obviously exact because $Z_nM,
B_nM\in \C$. Therefore, $\mathcal{L}^{\perp}\supseteq
(\dgF)^{\perp}=\barC$.

Let us now show the converse. Let $N\in \mathcal{L}^{\perp}$. We
must check that $N$ is exact and $Z_nN\in \C$. We first show that
$N$ is exact. This is equivalent to prove that any morphism
$S^n(G)\to N$ can be extended to a morphism $D^n(G)\to N$, for every
$n\in \Z$. But this fact follows from the short exact sequence
$$0\to S^n(G)\to D^n(G) \to S^{n-1}(G) \to 0 $$
since $\Ext^1(S^{n-1}(G),N)=0$.

Let us now check that $Z_nN\in \C$. We only need to show that
$\Ext^1_{\MMod}(F,Z_n N)=0$, for any $n\in \Z$ and any $F\in
S_{\aleph}$, since $S_{\aleph}$ cogenerates the cotorsion pair
$(\F,\C)$. By Lemma \ref{mist} we have a monomorphism of abelian
groups
$$0\to \Ext^1_{\MMod}(F,Z_n N)\to
\Ext^1_{\Ch(\MMod)}(S^n(F),N)$$ and, as we are assuming that
$\Ext^1_{\Ch(\MMod)}(S^n(F),N)=0$, we get that $Z_n N\in \C$. Hence
the cotorsion pair $(\dgF,\barC)$ is complete by
Theorem~\ref{bSalce}. \end{proof}

Now, we will prove that $(\barF,\dgC)$ is complete. We will first
need the following lemma. If $L=(L_i)_{i\in\Z}$ is a
complex, we will denote by $|L|$ the cardinality of $\oplus_{i\in\Z}L_i$.

\begin{lemma} \label{fccoconj}
Let $G$ be an h-unitary flat generator of $\MMod$ and let us fix $\aleph\geq
{\rm{max}}\{|G|,\aleph_0\}$. For any complex $F\in \barF$ and any
element $x\in F^k$ ($k\in\Z$ arbitrary), there exists a subcomplex
$L$ of $F$ such that $x\in L^k$, $|L|\leq \aleph$ and $L$, $F/L$ are
in $\F$.
\end{lemma}
\begin{proof} We may assume without loss of generality that $k=0$
and $x\in F^0$. By Lemma \ref{genpl}, there exists a morphism
$h_x:G\to F^0$ and an element $z_x\in G$ such that $h(z_x)=x$.
Consider then the exact complex
$$(S1)\hskip 1cm \cdots \rightarrow
A^{-2}_1\stackrel{\delta^{-2}}{\rightarrow} A^{-1}_1
\stackrel{\delta^{-1}}{\rightarrow} h_x(G)
\stackrel{\delta^{0}}{\rightarrow} \delta^0(h_x(G))
\stackrel{\delta^1}{\rightarrow} 0$$ where $A^{-i}_1$ is the submodule
of $F^{-i}$ constructed as follows. We know that $|h_x(G)|\leq
\aleph$, since $|G|\leq\aleph$. So we may find an $A^{-1}_1\leq
F^{-1}$ such that $|A^{-1}_1|\leq\aleph$ and $\delta^{-1}
(A^{-1}_1)=\ker (\delta^0|_{h_x(G)})$. Let us now choose an
$A^{-2}_1\leq F^{-2}$ with $|A^{-2}_1|\leq\aleph$ and
$\delta^{-2}(A^{-2}_1)=\ker (\delta^{-1} |_{A^{-1}_1})$. We
recursively repeat this argument for constructing all the $A^{-n}$.

On the other hand, $\ker(\delta^0|_{h_x(G)})\leq\ker \delta^0$ and
we know that $|\ker(\delta^0|_{h_x(G)})|\leq\aleph $. So by Lemma
\ref{tip}, $\ker (\delta^0|_{h_x(G)})$ can be embedded in a pure and
h-unitary  submodule $S^0_2$ of $\ker \delta^0$ in such a way that
$|S^0_2|\leq\aleph$. Consider the exact complex
$$(S2)\hskip 1cm \cdots \rightarrow
A^{-2}_2\stackrel{\delta^{-2}}{\rightarrow} A^{-1}_2
\stackrel{\delta^{-1}}{\rightarrow} h_x(G) +S^0_2
\stackrel{\delta^{0}}{\rightarrow} \delta^0(h_x(G))
\stackrel{\delta^1}{\rightarrow} 0$$ where the $A^{-i}_2$ $\!$'s are
constructed as above. Clearly $\ker (\delta^0|_{h_x(G)+S^0_2})
=S^0_2$, which is a pure and h-unitary  submodule of $\ker
\delta^0$. Moreover, $|h_x(G)+S^0_2|\leq\aleph +\aleph =\aleph$.

As $\delta^0(h_x(G))\subseteq \ker \delta^1$, we can embed
$\delta^0(h_x(G))$ in a pure and h-unitary  submodule $S^1_3$ of
$\ker\delta^1$ with $|S^1_3|\leq \aleph$, since
$|\delta^0(h_x(G))|\leq\aleph$. Let us consider the exact complex
$$(S3)\hskip 1cm \cdots \rightarrow
A^{-2}_3\stackrel{\delta^{-2}}{\rightarrow} A^{-1}_3
\stackrel{\delta^{-1}}{\rightarrow} A^0_3
\stackrel{\delta^{0}}{\rightarrow} S^1_3
\stackrel{\delta^1}{\rightarrow} 0.$$ We have again that $\ker
(\delta^1|_{S^1_3}) =S^1_3$, which is a pure and h-unitary
submodule of $\ker\delta^1$.

We turn over and find a pure and h-unitary  submodule
$S^0_4\leq\ker\delta^0$  with $|S^0_4|\leq\aleph$ and
$S^0_4\supseteq\ker (\delta^0|_{A^0_3})$. And then construct
$A^{-i}_4\leq F^{-i}$ ($|A^{-i}_4|\leq\aleph\ \forall i$) such that
$$(S4)\hskip 1cm \cdots \rightarrow
A^{-2}_4\stackrel{\delta^{-2}}{\rightarrow} A^{-1}_4
\stackrel{\delta^{-1}}{\rightarrow} A^0_3 +S^0_4
\stackrel{\delta^{0}}{\rightarrow} S^1_3
\stackrel{\delta^1}{\rightarrow} 0$$ is exact. Once more, $\ker
(\delta^0|_{A^0_3+S^0_4})= S^0_4\leq\ker\delta^0$ is a pure and
h-unitary submodule. Then we can find a pure and h-unitary submodule
$S^{-1}_5\leq\ker\delta^{-1}$ with $|S^{-1}_5|\leq\aleph$ and $\ker
(\delta^{-1}|_{A^{-1}_4})\subseteq S^{-1}_5$. Let us now consider the exact
complex $$(S5)\hskip 1cm \cdots \rightarrow
A^{-2}_5\stackrel{\delta^{-2}}{\rightarrow} A^{-1}_4+ S^{-1}_5
\stackrel{\delta^{-1}}{\rightarrow} A^0_3 +S^0_4
\stackrel{\delta^{0}}{\rightarrow} S^1_3
\stackrel{\delta^1}{\rightarrow} 0,$$ in which
$\ker(\delta^{-1}|_{A^{-1}_4+S^{-1}_5})
=S^{-1}_5\leq\ker\delta^{-1}$ is pure and h-unitary .

Our next step will be to find a pure and h-unitary  submodule
$S^{-2}_6\leq\ker\delta^{-2}$ such that $|S^{-2}_6|\leq\aleph$ and
$\ker(\delta^{-2}|_{A^{-2}_5})\subseteq S^{-2}_6$. And then
consider the exact complex $$(S6)\hskip 1cm \cdots \rightarrow
A^{-3}_6\stackrel{\delta^{-3}}{\rightarrow}
A^{-2}_5+S^{-2}_6\stackrel{\delta^{-2}}{\rightarrow} A^{-1}_4+
S^{-1}_5 \stackrel{\delta^{-1}}{\rightarrow} A^0_3 +S^0_4
\stackrel{\delta^{0}}{\rightarrow} S^1_3
\stackrel{\delta^1}{\rightarrow} 0$$ in which
$\ker(\delta^{-2}|_{A^{-2}_5+S^{-2}_6})=S^{-2}_6 \subseteq
\ker\delta^{-2}$ is a pure and h-unitary  submodule.

Finally, we can prove by induction that, for any $n\geq 4$, we may
construct an exact complex $$(Sn)\cdots
\stackrel{\delta^{-n+2}}{\rightarrow} A^{-n+3}_n
\stackrel{\delta^{-n+3}}{\rightarrow} T^{-n+4}_n
\stackrel{\delta^{-n+4}}{\rightarrow} T^{-n+5}_n \rightarrow \cdots
\stackrel{\delta^{-1}}{\rightarrow}
T^0_n\stackrel{\delta^0}{\rightarrow} T^1_n
\stackrel{\delta^1}{\rightarrow} 0$$ such that $\ker
(\delta^{-n+j}|_{T^{-n+j}_n})$ is a pure and h-unitary  submodule of
$\ker\delta^{-n+j}$ $\forall j\geq 4$ and all the terms have
cardinality bounded by $\aleph$.

Let $L=\underset{n\in\N}{\underrightarrow{\text{lim}}} (Sn)$. It is
straightforward to check that $L$ is an exact complex (since it is a
direct limit of exact complexes). And that $\ker(\delta^i|_{L^i})$
is a pure and h-unitary  submodule of $\ker\delta^i$ $\forall i\leq
1$. Furthermore $|L^i|\leq\aleph_0\cdot\aleph= \aleph$ for any
$i\leq 1$ and thus, $|L|\leq\aleph$.

This complex
$$L=\cdots \rightarrow L^i\stackrel{\delta^i}{\rightarrow} L^{i+1}
\stackrel{\delta^{i+1}}{\rightarrow} \cdots
\stackrel{\delta^{-1}}{\rightarrow}
L^0\stackrel{\delta^0}{\rightarrow} L^1
\stackrel{\delta^1}{\rightarrow} 0 \stackrel{\delta^2}{\rightarrow}
0\cdots,$$ is a subcomplex of $F$ which satisfies that $x\in L^0$
and that $\ker(\delta^i|_{L^i})$ is a pure and h-unitary  submodule
of $\ker\delta^i$ for any $i\in\Z$. In particular,
$\ker(\delta^i|_{L^i})$ is an h-unitary flat module for any $i\in\Z$
(since so is $\ker\delta^i$). Therefore, the complex $L$ is an
h-unitary flat subcomplex of $F$ satisfying that $|L|\leq\aleph$.

To finish the proof we only need to check that
$F/L=(F^i/L^i,\overline{\delta}^i)$ is an h-unitary flat complex. An
easy computation shows that
$\ker\overline{\delta}^i=\ker(\delta^i)/\ker(\delta^i|_{L^i})$. But,
by construction, $\ker(\delta^i|_{L^i})$ is a pure and h-unitary
submodule of $\ker\delta^i$, $\forall i\in\Z$. So
$\ker(\delta^i)/\ker(\delta^i|_{L^i})$ is h-unitary  and flat for
any $i\in \Z$. Finally, $F/L$ is exact since both $F$ and $L$ are
exact. Thus, $F/L$ is an h-unitary  and flat complex.
\end{proof}

\begin{theorem}
The pair $(\barF,\dgC)$ is complete.
\end{theorem}
\begin{proof} Let $F$ be an h-unitary  and flat complex and choose an $x\in F^i$.
We can find, by the above result, a subcomplex $L\in \barF$ of $F$
such that $x\in L^i$, $|L|\leq \aleph$ and the quotient complex
$F/L$ is h-unitary  and flat. Then we can mimic the proof of Theorem
\ref{firmcomp} for showing that the cotorsion pair $(\barF,\dgC)$ is
complete.
\end{proof}

Let us denote by
$\J:\Firm\hookrightarrow R\Mod$ the embedding functor. As $\Firm$
has a generator (see Lemma \ref{genpl}) and $\J$ preserves small colimits, the Special
Adjoint Functor Theorem (see e.g. \cite{freyd}) ensures the
existence of a right adjoint functor $\D:R\Mod\to \Firm$ of $\J$.
Thus, $\Firm$ is a coreflective subcategory of $R\Mod$ and limits in
$\Firm$ are computed by applying the functor $\D$ to the usual
limits in $R\Mod$ (see e.g. \cite{Adamek}).

\begin{lemma}\label{et1}
$\dgF\cap \mathcal{E}=\barF$ and $\dgC\cap \mathcal{E}=\barC$ in
$\Ch(\MMod)$.
\end{lemma}
\begin{proof} Let us choose $X\in\barF$, and $Y\in \barC$.
Then any map from $X$ to $Y$ is homotopic to zero and so $Hom(X,Y)$
is exact. Therefore $X\in \dgF\cap \mathcal{E}$. Conversely, suppose
that $X\in\dgF$ is exact. We have to show that the h-unitary  cycle
module $Z_nX$ is flat, $\forall n\in \Z$. Let $C$ be any cotorsion
complex in $\Ch(R)$. As the embedding functor
$\J:\Ch(\Firm)\rightarrow \Ch(R)$ is exact (and thus, $\D(C)\in
\barC$), we deduce that its right adjoint functor $\D$ preserves
cotorsion complexes and therefore,  $Hom_{\Ch(R)}(\J(X),C)\cong
Hom_{\Ch(\Firm)}(X,\D(C))$ is exact. Hence, $\J(X)$ is a dg-flat and
exact complex in $R\Mod$. By the arguments used in the proof of
\cite[Lemma 4.4.8]{jr}, we deduce that $\J(X)$ is a flat complex and
so, $Z_n\J(X)=Z_nX$ is flat. Therefore $X\in \barF$.

Let us now check that $\dgC\cap \mathcal{E}=\barC$. As before, we
have that $\barC\subseteq \dgC\cap \mathcal{E}$. Conversely, assume
that $Y\in \dgC$ is exact. We must check that $Y\in \barC$. By
Theorem~\ref{comp1}, there exists a short exact sequence
$$0\to Y\to C\to H\to 0 $$ with $C\in \barC$ and $H\in \dgF$. Let us note that
$H$ is exact, since so are $Y$ and $C$. And $H\in \barF$ because
$H\in \dgF$. But then the above short exact sequence splits, by
Lemma \ref{cotpa}. So $Y$ is a direct summand of $C$. And this means
that $Y\in \barC$, since $\barC$ is closed under direct summands.
 \end{proof}


\section{A Unitless Monoidal Structure in $\Ch(\Firm)$}

We devote this section to show that the Flat Model Structure defined in last section on $\Ch(\MMod)$ is compatible
 with the tensor
product in $\Ch(\Firm)$ inherited from the tensor product in
$\Firm$. Note that we may assume that any left $R$-module is also a right $R$-module by \cite{loday}.

Let us note that $\Firm$ with the induced tensor product of $R\Mod$
has all the structure of a symmetric monoidal category apart from
the unit object. This is known in the literature as a {\it unitless monoidal categoy} (see \cite{hines}). Moreover we can get a closed
structure in $\Firm$ by applying the functor $\mathbf{D}$ after the
usual internal Hom functor of $R\Mod$. So $\Firm$ is a unitless
closed symmetric monoidal category. This imposes a canonical
unitless closed symmetric monoidal structure in $\Ch(\Firm)$.

In order to prove that the model structure defined in $\Ch(\MMod)$ in last section  is compatible with this tensor product of $\Ch(\Firm)$, we need to prove the following
technical Lemma. Recall that the class $\Lt$ appearing on it was
introduced in the proof of Theorem~\ref{comp1}.

\begin{lemma}\label{menc}
Let $M\in \dgF$. Then $M$ is a direct summand of an $\Lt$-filtered
dg-h-unitary flat complex.
\end{lemma}

\begin{proof}
We know that $(\dgF,\barC)$ is a cotorsion pair in $\Ch(\MMod)$ cogenerated by
$\Lt$, by Theorem~\ref{comp1}. Thus, $\dgF={^\perp}(\Lt^\perp)$ and
$\barC=\Lt^\perp$. By \cite{EkTrl}, for every complex of h-unitary
modules $K$, there exists a short exact sequence
$$0\to K\to Y\to Z\to 0$$ with $Y\in \Lt^{\perp}=\barC$ and $Z$, an $\Lt$-filtered
complex (in particular, $Z\in \dgF$).

We are going to finish the proof by adapting
\cite[Lemmas 2.2 and 2.3]{Salce}. Given any $M\dgF$,
there exists a short exact sequence
$$0\to K\to \oplus_{n\in \Z}S^n(G)\to M\to 0 $$ where $G$ is an h-unitary flat generator of $\MMod$.
Now let
$$0\to K\to Y\to Z\to 0$$ be any short exact sequence  with $Y\in \barC$ and $Z$,
$\Lt$-filtered. Let us construct the pushout \DIAGV{60} {} \n{}
\n{0}\n{}\n{0}\nn {}\n{}\n{\sar}\n{}\n{\sar}\nn{0}
\n{\ear}\n{K}\n{\ear}\n{\oplus_{n\in
\Z}S^n(G)}\n{\ear}\n{M}\n{\ear}\n{0}\nn
{}\n{}\n{\sar}\n{}\n{\sar}\n{}\n{\seql}\nn
{0}\n{\ear}\n{Y}\n{\ear}\n{Q}\n{\ear}\n{M}\n{\ear}\n{0}\nn
{}\n{}\n{\sar}\n{}\n{\sar}\nn {}\n{}\n{Z}\n{\eeql}\n{Z}\nn
{}\n{}\n{\sar}\n{}\n{\sar}\nn {}\n{}\n{0}\n{}\n{0}\diag
As $\oplus_{n\in \Z}S^n(G)$ is trivially $\Lt$-filtered and $Z$ is
$\Lt$-filtered, we see that $Q$ is also $\Lt$-filtered. And clearly
$Y\in \barC=\Lt^{\perp}$. Hence, as $M\in
\dgF={^\perp}(\Lt^{\perp})$, we get that the short exact sequence
$$0\to Y\to Q\to M\to 0$$
splits. And thus, $M$ is a direct summand of the $\Lt$-filtered
complex $Q$. \end{proof}

We can now state the main result of this section.
\begin{theorem}\label{escom}
The model structure induced on $\Ch(\MMod)$ by the flat cotorsion
pair $(\F,\C)$ is compatible with the symmetric closed unitless monoidal
structure of $\Ch(\Firm)$ induced by the usual tensor product
of chain complexes.
\end{theorem}
\begin{proof} We are going to use \cite[Theorem 7.2]{hovey2}, so we will need to check
that
\begin{enumerate}
\item Every cofibration is pure.
\item If $X,Y\in \dgF$, then $X\otimes Y\in \dgF$.
\item If $X,Y\in \dgF$ and one of the them belongs to $\barF$, then $X\otimes Y\in \barF$.
\end{enumerate}
As any cofibration is a monomorphism with cokernel in $\dgF$, it
follows that it is a short exact sequence with h-unitary flat
cokernel. Thus, it is a pure exact sequence. This shows that $(1)$
holds.

Let us choose $X,Y\in \dgF$. By Lemma~\ref{menc}, we only need to
check condition $(2)$ for complexes of the form $S^n(F)$, with $F$ a
h-unitary flat module. But in this case, $S^n(F)\otimes S^m(F')\cong
S^{n+m}(F\otimes F')$, for any other h-unitary flat module $F'$. And, as the tensor product of two h-unitary
flat modules is clearly h-unitary  and flat, we conclude that
$X\otimes Y$ is a direct summand of the direct limit of the directed
system of complexes $S^t(G\otimes G')$ in $\dgF$ and hence, it is
also complex in $\dgF$.

Let us now check condition $(3)$. Suppose that $X\in \barF$ and
$Y\in \dgF$. $X\otimes Y$ is exact, because $X$ is exact. And
$X\otimes Y\in \dgF$, by the previous remarks. Therefore $X\otimes
Y\in \barF$ by Lemma~\ref{et1}.
\end{proof}

\section{Morita Invariance}\label{moritinv}

We finish this paper by showing the invariance of the previous
constructions under the equivalences induced by the Morita contexts considered in \cite{Quillen2}. Namely, the definitions of $\MMod$ and $\Firm$ apparently depend on the fact that we fixed the embedding of $A$ in the extended ring $\tilde{A}$. But we are going to prove in this section that we get isomorphic definitions of $\MMod$ and $\Firm$ if we consider any other embedding of $A$ as a two-sided ideal of a unital ring $R$.

\begin{theorem}\label{esMi}
Let $A$ be a two-sided ideal of a unital ring $R$. Let us denote by $\MMod_R$,
$\F_R$ and $\C_R$ the categories of h-unitary $R$-modules, h-unitary
flat $R$-modules and h-unitary cotorsion $R$-modules respectively.
Then there exist canonical equivalences $\MMod_R\cong \MMod$,
$\F_R\cong \F$ and $\C_R\cong \C$.
\end{theorem}

\begin{proof}
Let $M\in \MMod_R$. By the comments of the previous sections, there
exists an $\F_R$-resolution of $M$

$$\cdots\to F_1\to F_0\to M\to 0.$$

On the other hand, it is clear that $A\otimes_{\tilde{A}}M\cong
A\otimes_R M$ for all $R$-module $M$ such that $AM=M$ (in
particular for all $M\in \MMod_R$). Hence
$\Tor_j^{\tilde{A}}(A,M)=0$, $\forall j\geq 1$ and so $M\in \MMod$.
We get then that, for a given $R$-module $M$, $M\in \MMod_R$ if and
only if $M\in \MMod$.

Conversely, given $M$ in $\MMod$, we have that $A\otimes_{\tilde{A}}M\cong
M$. So $A\otimes_{\tilde{A}}M$ has an $R$-module structure given
by $r(a\otimes m)=ra\otimes m$ and therefore, $M$ has a unique $R$-module
structure extending its $\tilde{A}$-module structure. By the
preceding paragraph it follows that $M\in \MMod_R$. Hence we have a
one-to-one correspondence between $h$-unitary module structures on any
abelian group for the pairs $(R,A)$ and $(\tilde{A},A)$.

Now let us choose $F\in \F$. We have that $AF=F$. So there is a canonical
isomorphism $M\otimes_{\tilde{A}}F\cong M\otimes_R F$ for all right
$R$-module $M$. It follows that $F\in \F_R$. Conversely, assume that $F'
\in \F_R$ and let $N\to M$ be an injection of right
$\tilde{A}$-modules. We get an exact sequence of right $R$-modules

$$\Tor_1^{\tilde{A}}(M/N,R)\to N\otimes_{\tilde{A}}R\to
M\otimes_{\tilde{A}}R.$$

But it is easy to check that
$\Tor_1^{\tilde{A}}(M/N,R)A=0$. So applying the exact functor
$-\otimes_{R}F'$ to the previous exact sequence we get that
$N\otimes_{\tilde{A}}F'\to M\otimes_{\tilde{A}}F'$ is injective and
hence, $F'\in \F$.

The analogous statement for the subcategories $\C_R$ and $\C$ is now
easy to prove by using the preceding paragraph and noting that an exact
sequence $0\to N\to M\to M/N\to 0$ in $\MMod$ is splitting if and only
if it splits in $\MMod_R$.
\end{proof}

\begin{corollary}
The previous h-unitary flat model structure on $\Ch(\MMod)$ is independent of the embedding of $A$ into a unital ring.
\end{corollary}

\noindent {\bf ACKNOWLEDGEMENT.}
The authors wish to thank Prof. Edgar Enochs for his comments
during the preparation of this work.

\end{document}